%% file: 00-bdcok-preprint.tex
\documentclass[a4paper]{article}
\usepackage{newclude}
\usepackage{graphicx} 
\usepackage{amsmath, amsthm, amsfonts, amssymb, amscd}
\usepackage{setspace}
\usepackage{enumitem}
\usepackage{latexsym}
\usepackage{verbatim}
\usepackage[colorlinks=true,hypertexnames=false,linkcolor=blue,citecolor=blue, backref=page]{hyperref}
\usepackage{mathrsfs}
\usepackage{color}
\usepackage[usenames,dvipsnames]{xcolor}
\usepackage{cite}
\usepackage[noblocks]{authblk}
\usepackage[utf8]{inputenc}

\pdfoutput=1
\graphicspath{{figs/}}

\advance \textheight by \topskip
\numberwithin{equation}{section}
\usepackage{fancyhdr}
\def\inn{{inn}}
\def\out{{out}}

\newtheorem{mtheorem}{Theorem}
\newtheorem{theorem}{Theorem}[section]
\newtheorem{lemma}[theorem]{Lemma}
\newtheorem{proposition}[theorem]{Proposition}
\newtheorem{corollary}[theorem]{Corollary}

\theoremstyle{definition}
\newtheorem{definition}[theorem]{Definition}
\newtheorem{remark}[theorem]{Remark}

\makeatletter
\renewenvironment{proof}[1][\proofname]{\par
\pushQED{\qed}%
\normalfont \topsep6\p@\@plus6\p@\relax
\trivlist
\item\relax
{\bfseries
#1\@addpunct{.}}\hspace\labelsep\ignorespaces
}{%
\popQED\endtrivlist\@endpefalse
}
\makeatother
\def\ackname{Acknowledgements}%
\newenvironment{ack}[1][\ackname]%
{\ifx#1\empty\else\subsection*{#1.}\fi\par}
{\par}
\usepackage{cancel}
\usepackage[normalem]{ulem}
\newcommand\redst{\bgroup\markoverwith{\textcolor{red}{\rule[0.5ex]{2pt}{0.7pt}}}\ULon}
\newcommand\bluest{\bgroup\markoverwith{\textcolor{blue}{\rule[0.5ex]{2pt}{0.7pt}}}\ULon}
\usepackage[textwidth=2.7cm, color=red!20]{todonotes}
\newcounter{nota}[section]

\newcommand{\ver}{} 

\title{Hyperbolicity and Abundance of Elliptical Islands in Annular Billiards}
\author{R B Batista}
\affil{\small Departamento de Matemática, ICE, UFJF, Brasil (reginaldo.braz@ufjf.br)}
\author{M J Dias Carneiro}
\author{S Oliffson Kamphorst}
\affil{\small Departamento de Matemática, ICEX, UFMG, Brasil (carneiro@mat.ufmg.br, syok@ufmg.br)}
\date{}
\begin{document}
\maketitle

\include*{0-abs}
\include*{1-intro}

\include*{2-preliminares}
\include*{3-cones}
\include*{4-pontosnormais}
\include*{5-hyperbolic}

\include*{6-tangencias}
\begin{ack}
We want to thank S Pinto-de-Carvalho for her contributions. 
This work originated from the thesis of RB Batista \cite{tese} with the funding of
Coordena\c c\~ao de Aperfei\c coamento de Pessoal de N\'\i vel Superior (CAPES) and 
Conselho Nacional de Desenvolvimento Cient\'\i fico e Tecnol\'ogico (CNPq),  Brasil. 
\end{ack}
\bibliographystyle{abbrv}
\bibliography{ab}

 \end{document}

%% file: 0-abs.tex
\begin{abstract}
We study the billiard dynamics in annular tables between two excentric circles. 
As the center and the radius of the inner circle change, a two parameters map is defined by the first return of trajectories to the obstacle. 
We obtain an increasing family of hyperbolic sets, in the sense of the Hausdorff distance, 
as the radius goes to zero and the center of the obstacle approximates the outer boundary. 
The dynamics on each of these sets is conjugate to a shift with an increasing number of symbols. 
We also show that for many parameters the system presents quadratic homoclinic tangencies whose bifurcation gives rise to elliptical islands
(Conservative Newhouse Phenomenon).
Thus, for many parameters we obtain the coexistence of a "large" hyperbolic set with many elliptical islands. 
\end{abstract}

%% file: 1-intro.tex
\section{Introduction}

The billiard problem consists of the description of the free motion of a particle 
inside a region of the plane called the {\em table}. The particle moves  
in straight lines with constant unitary speed between the boundaries and undergoes elastic collisions at the impacts. Conservation of energy and momentum implies the reflection 
law for the collisions with the boundaries. 
The two dimensional dynamics is given by the  {\em billiard map}, assigning a collision to the next one.
The dynamical properties of a billiard, which are deeply related to its shape, range from integrability to ergodicity.

The mathematical billiards were introduced by Birkhoff \cite{birkhoff}  who showed that  the motion on elliptical tables is integrable, he also conjectured
that these are the only convex billiards with this property \cite{conjectura}. 
Birkhoff billiards, as strictly convex billiards are now called, in general exhibit invariant curves and elliptical islands coexisting with regions of hyperbolic behavior.
A full description of the dynamics in generic convex billiards is still a challenge.
On the other hand, Sinai \cite{sinai} used dispersing billiards to investigate the micro dynamics of the ideal gas and address the Boltzmann's ergodic hypothesis. The so-called Sinai's billiard, which is a classical example of a dispersive billiard, was proved to be chaotic.
It is now known that if 
all the components of the boundary of the table are concave the dynamics is fully chaotic,
i.e with positive Lyapunov exponent a.e., and in general can be shown to be ergodic.
So, the common idea it that billiards with concave/dispersing components are associated to  hyperbolic/random/chaotic behavior,
while convex/focusing components frequently imply some non chaotic/elliptical  behavior with regions of stability.

In this work we study billiards in annular tables as introduced by Saito et al. in \cite{saito}.
An annular table, is {\ver a} closed planar region $Q_{\delta,r} \subset {\mathbb R}^2$  bounded by an unitary circle $\gamma$ centered at the origin, and an inner circle  $\alpha$ of radius $r$ 
centered at a point $p_{\delta}$ at a distance $\delta$ from the origin.
We  call the exterior unitary circle the {\em exterior boundary} and the inner circle is called the {\em obstacle}. The distance $\delta$ between the centers is called the {\em eccentricity}.  
The set of parameters is the triangular region 
$ \Omega=\{(\delta,r): 0\leq \delta <1\ \,  \hbox{ and }  0< {\ver r + \delta} <1\}$.
The corresponding two parameter family of billiard maps is denoted by $T_{\delta,r}$.
As the collisions with the inner circle carry the interesting part of the dynamics, 
it is meaningful to describe the dynamics through the first return to the obstacle map, denoted by $G_{\delta,r}$. 
Our results are stated for this first return map and correspondingly we refer, in this introduction, to the set of collisions with the obstacle as the phase space.
\begin{figure}[h]
{\includegraphics[width=0.3\hsize]{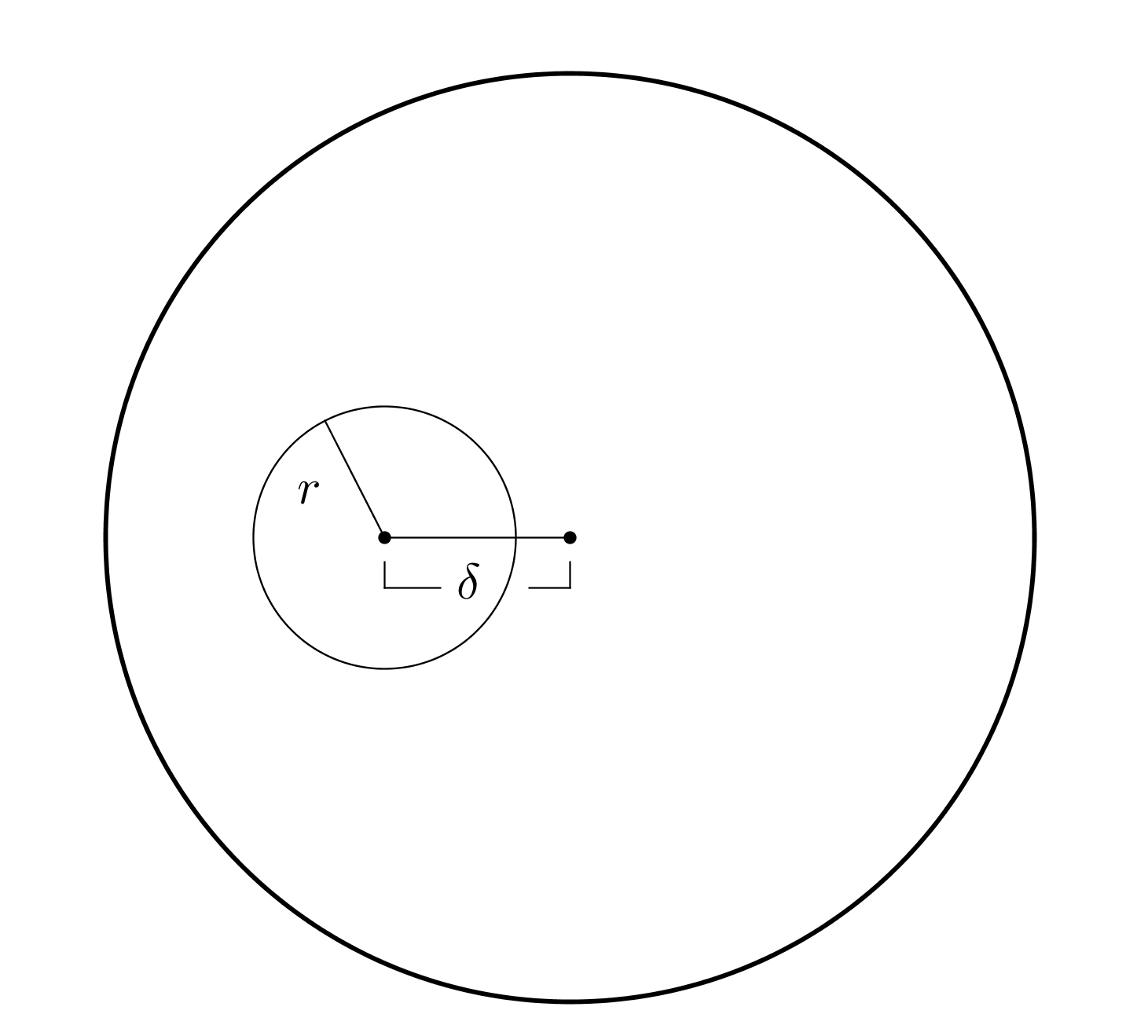}}\hfill{\includegraphics[width=0.65 \hsize]{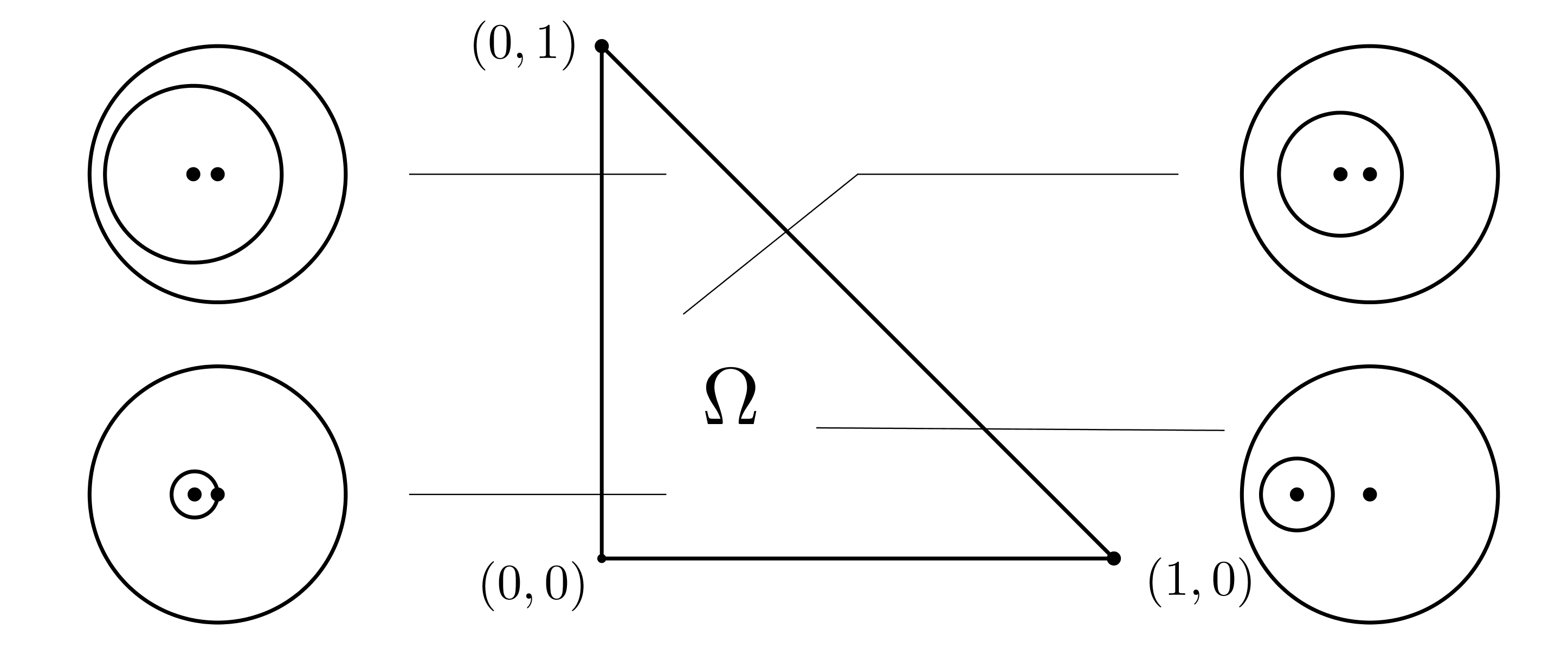}}
\caption{The annular table and the parameter space $\delta \times r$}
\end{figure}

\begin{figure}[h]
{
\includegraphics[width=0.19\hsize]{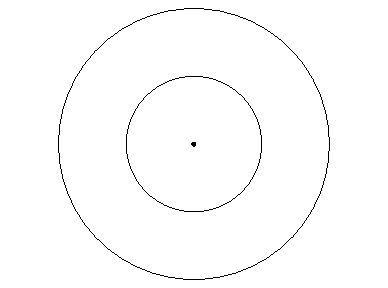}\hfill
\includegraphics[width=0.19\hsize]{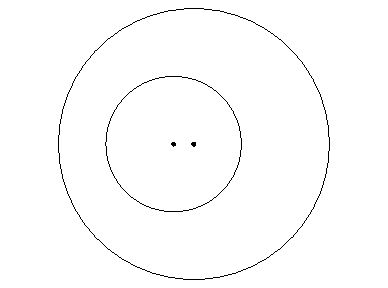}\hfill
\includegraphics[width=0.19\hsize]{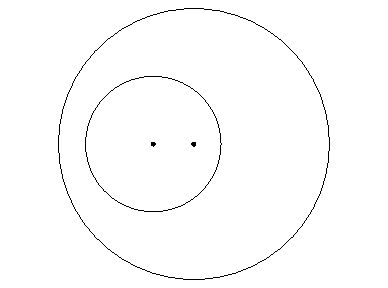}\hfill
\includegraphics[width=0.19\hsize]{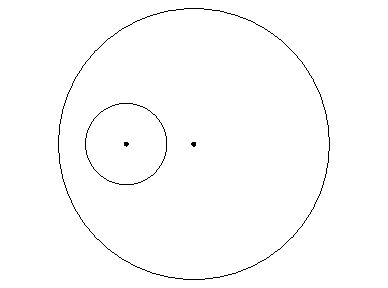}\hfill
\includegraphics[width=0.19\hsize]{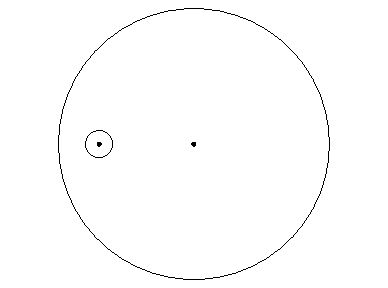}
}
\vskip 0.5cm
{
\includegraphics[width=0.19\hsize]{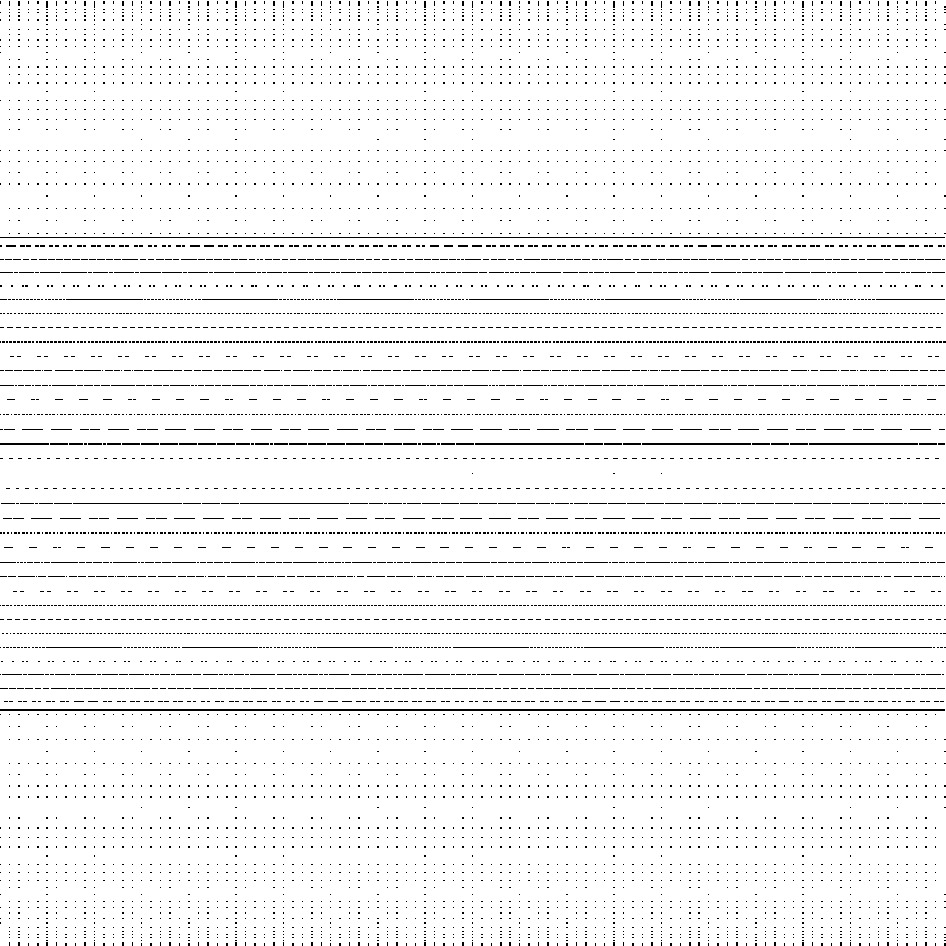}\hfill
\includegraphics[width=0.19\hsize]{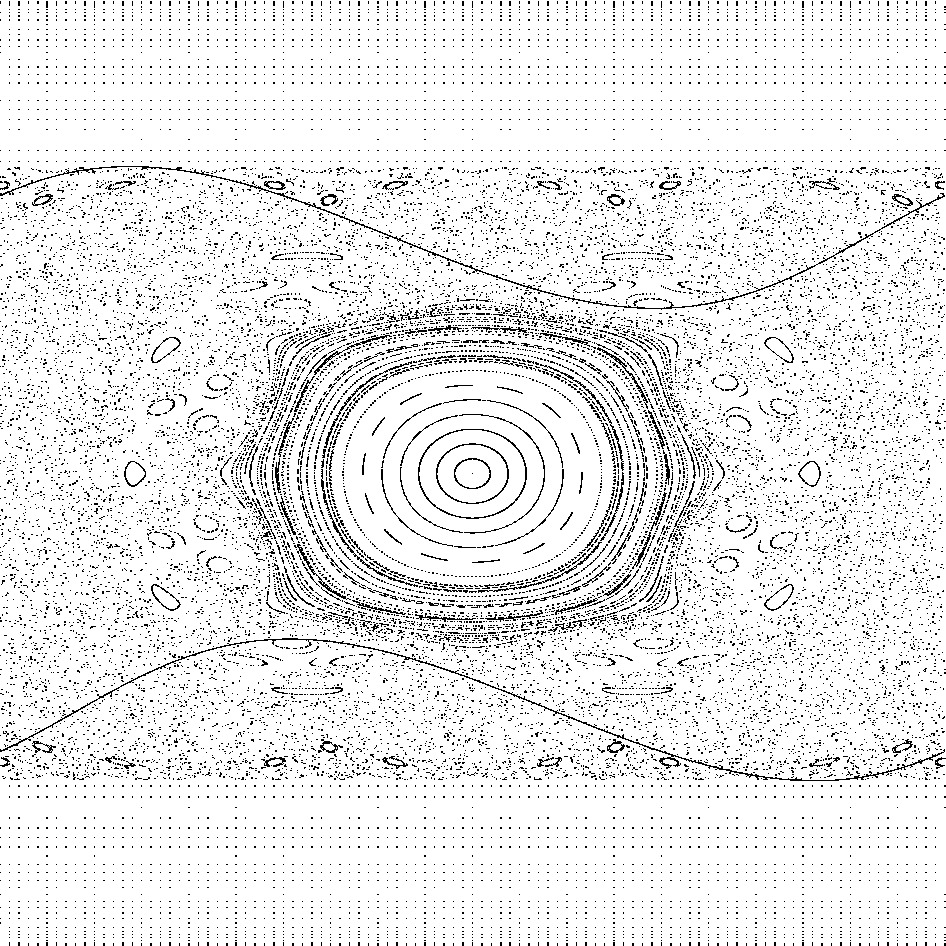}\hfill
\includegraphics[width=0.19\hsize]{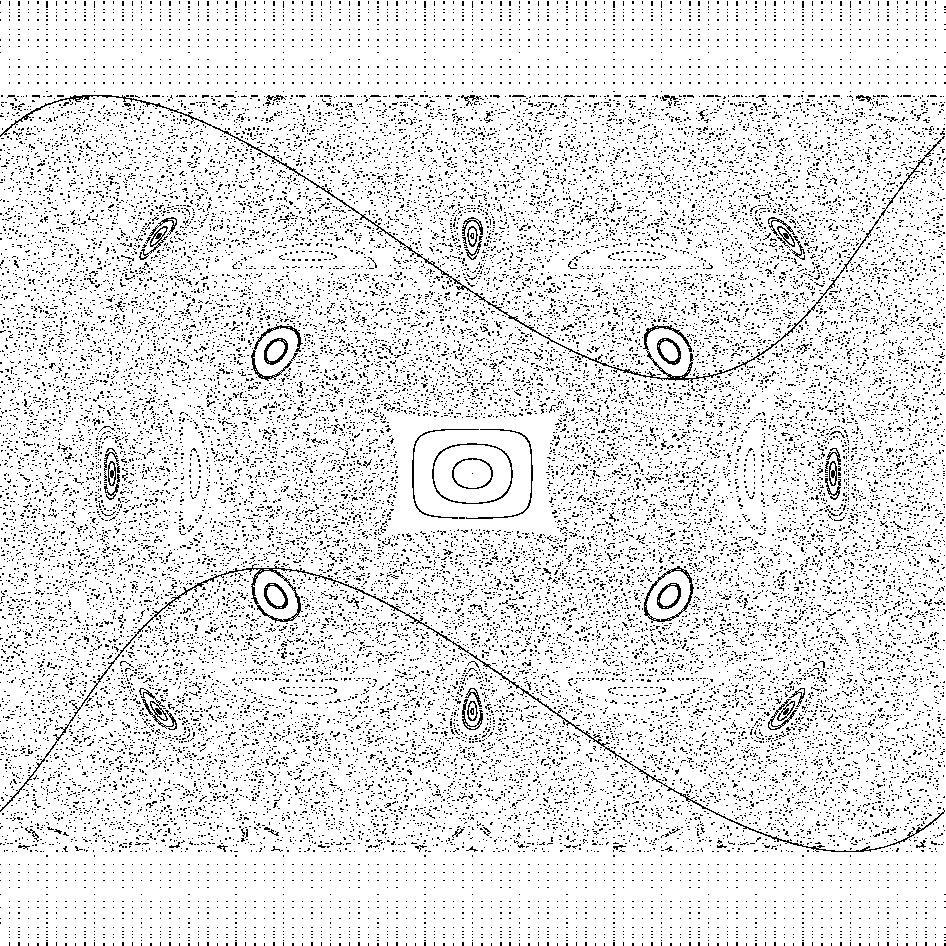}\hfill
\includegraphics[width=0.19\hsize]{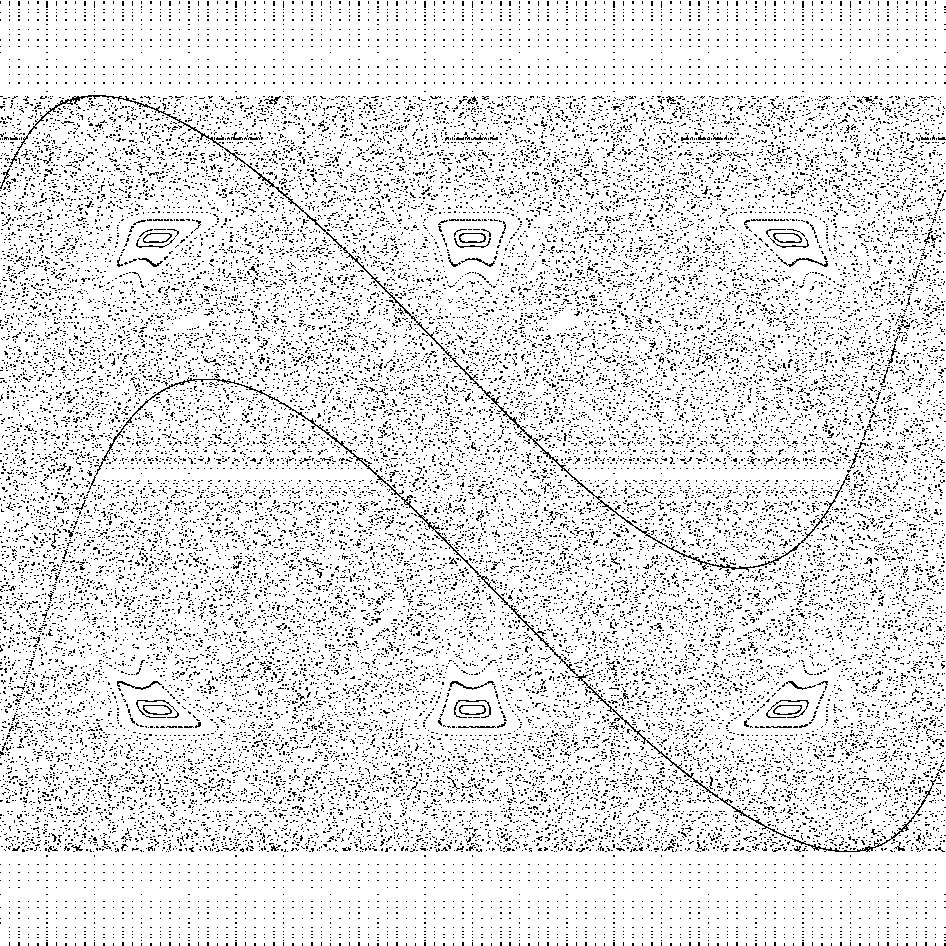}\hfill
\includegraphics[width=0.19\hsize]{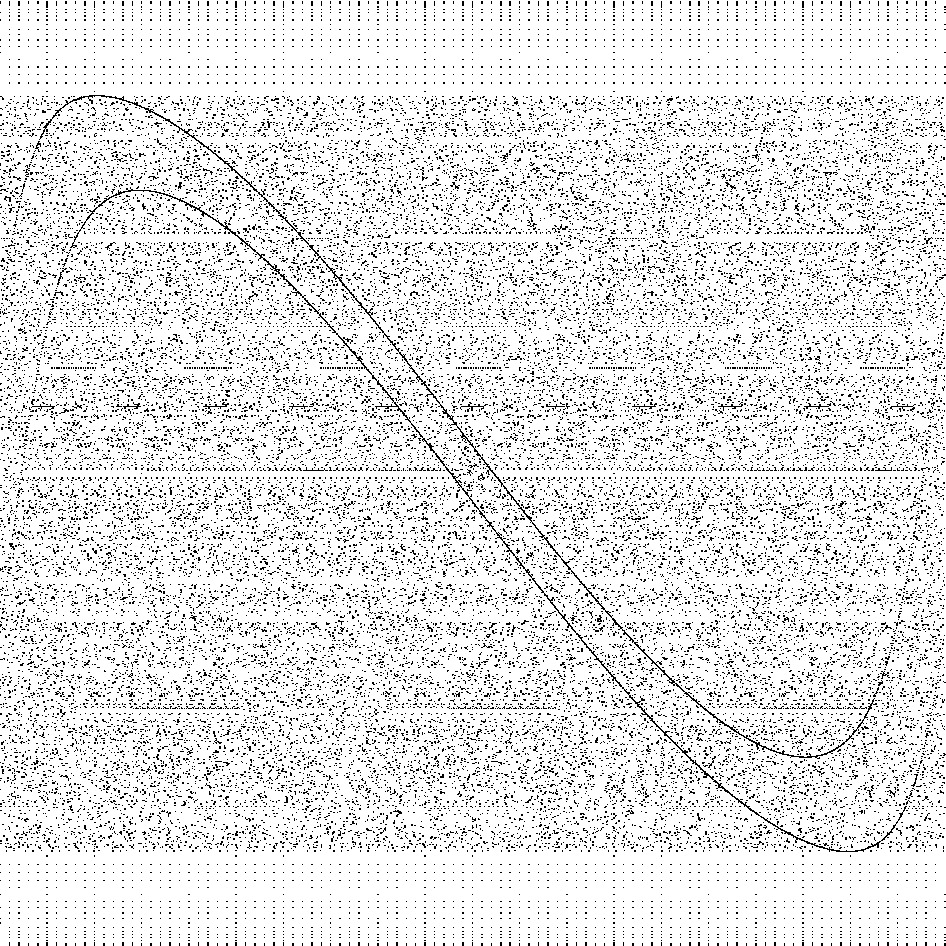}
}\vskip 0.5cm
{
\includegraphics[width=0.19\hsize]{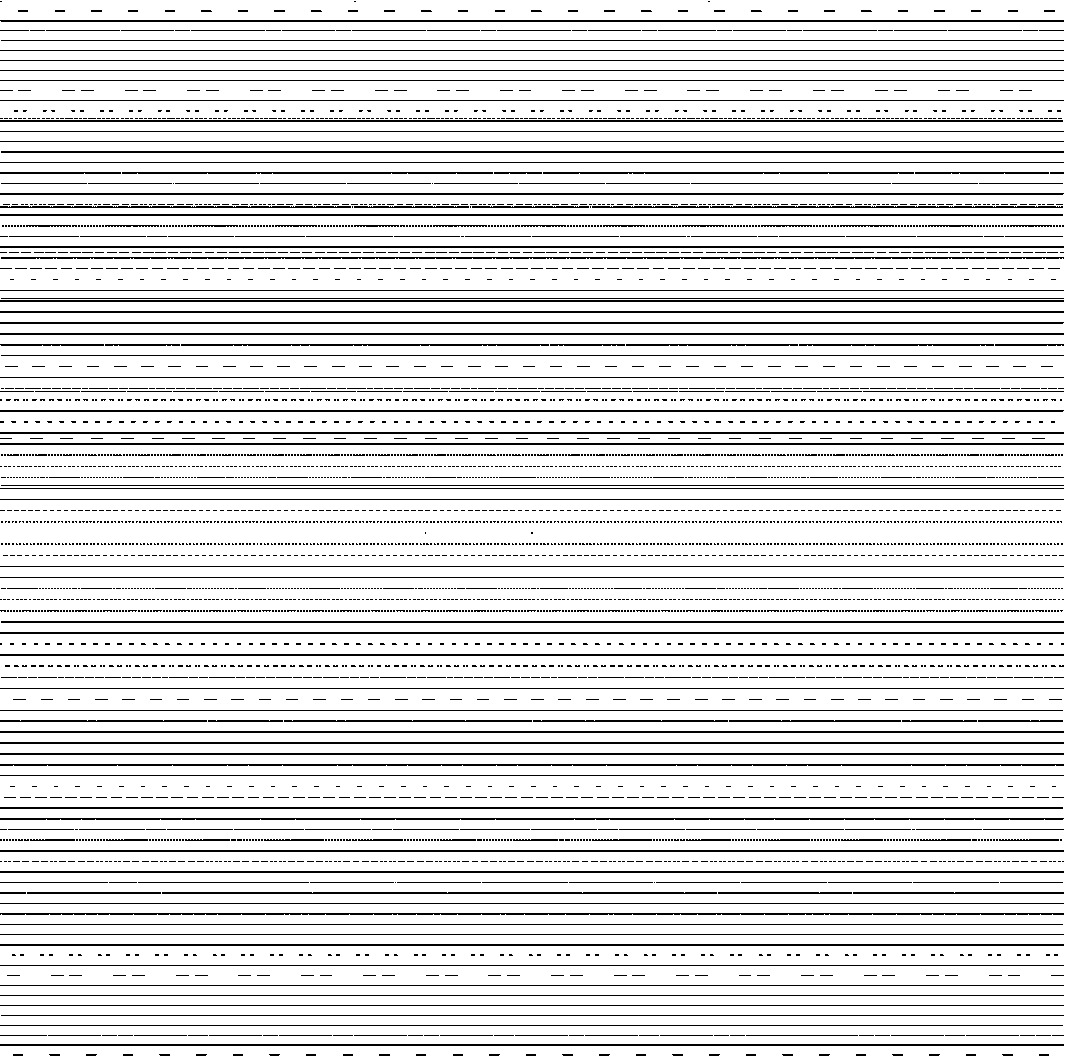}\hfill
\includegraphics[width=0.19\hsize]{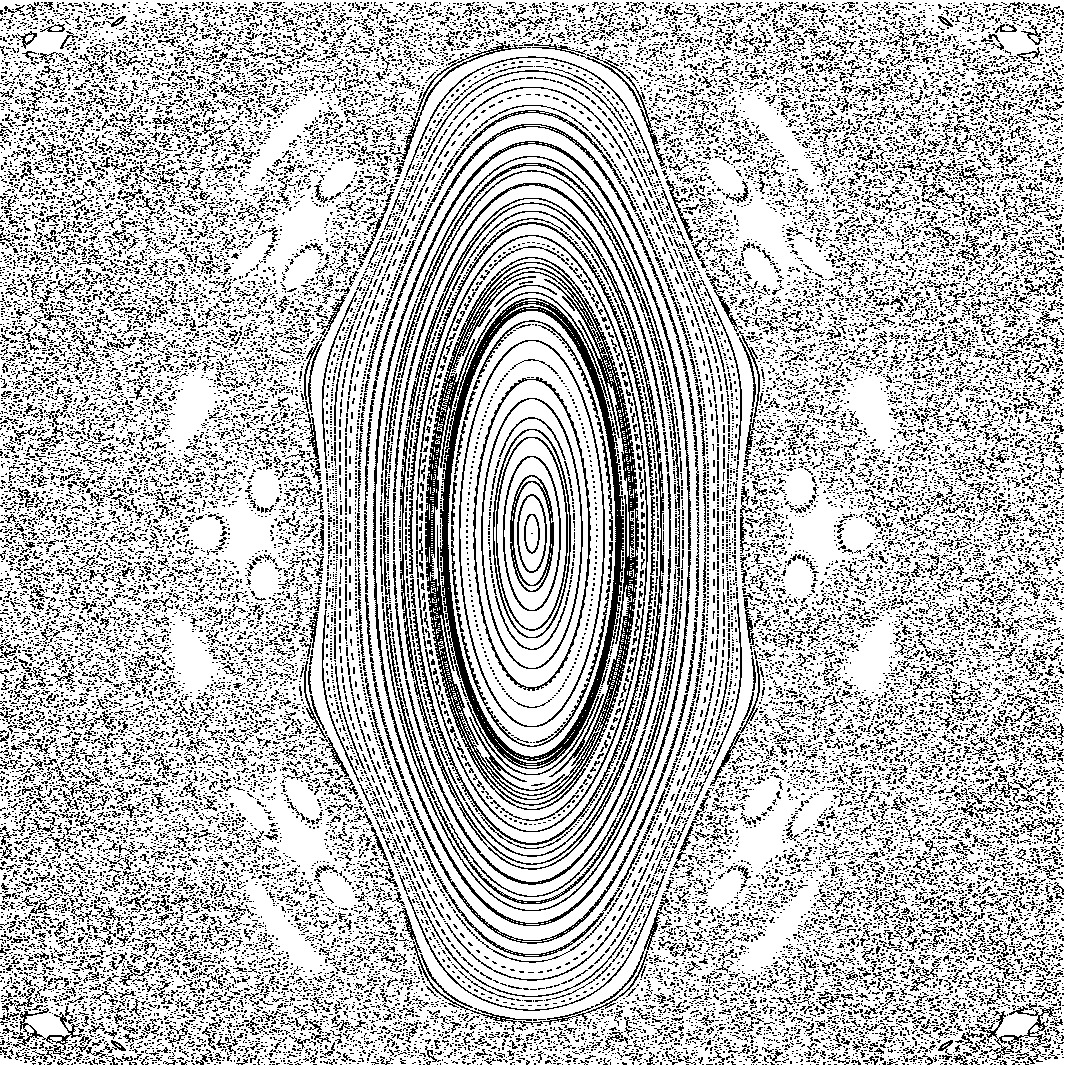}\hfill
\includegraphics[width=0.19\hsize]{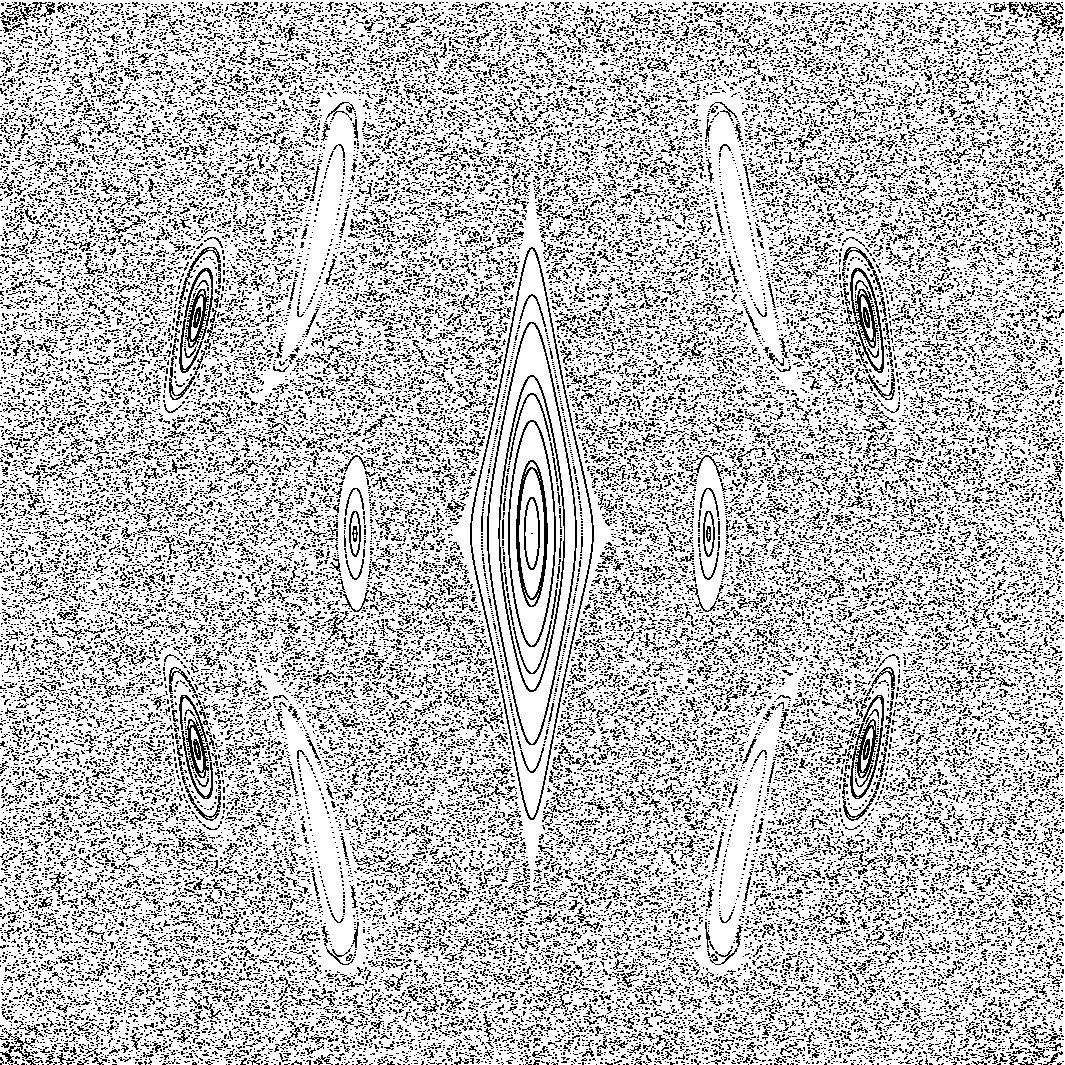}\hfill
\includegraphics[width=0.19\hsize]{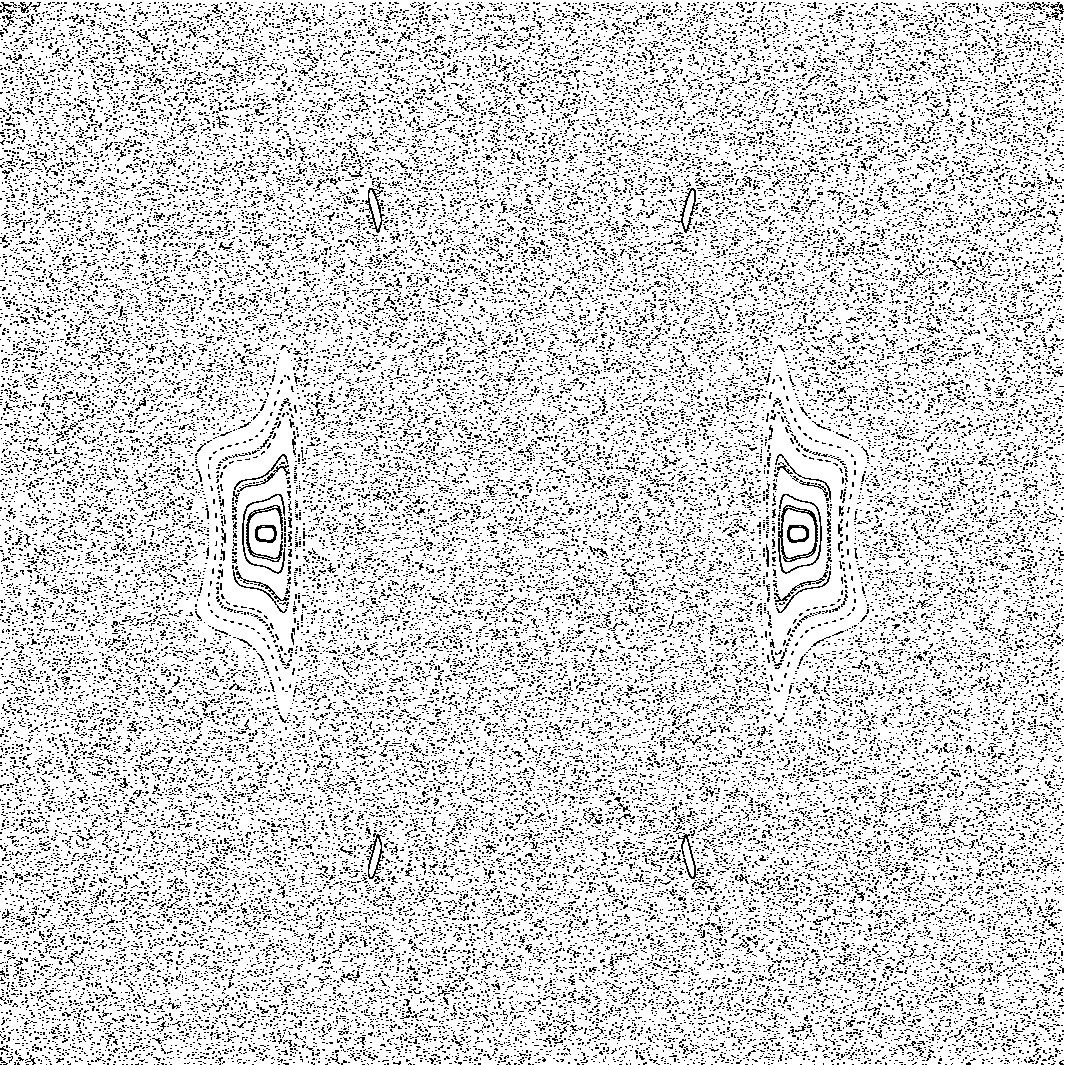}\hfill
\includegraphics[width=0.19\hsize]{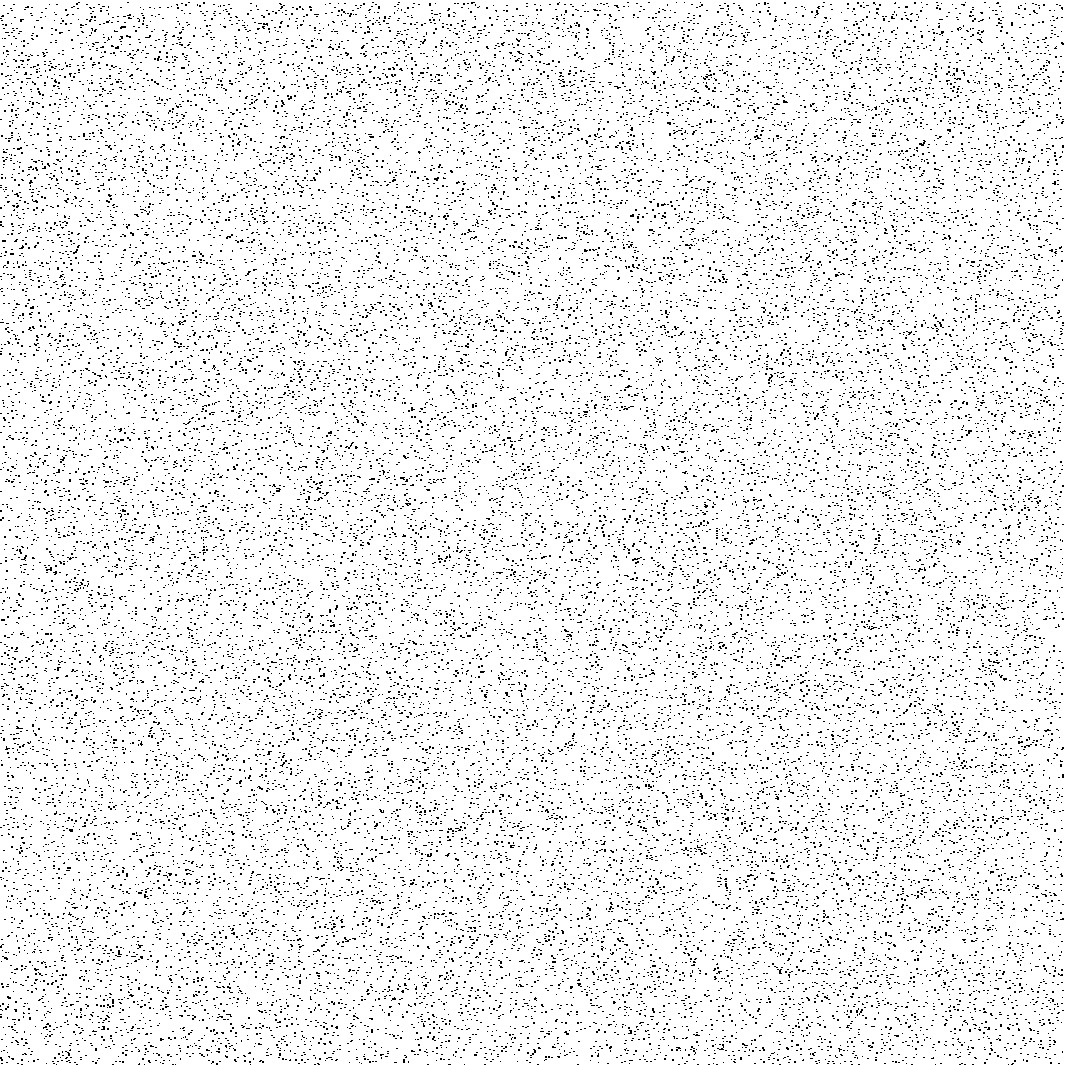}
}
\caption{Change of the dynamics with the parameters as viewed from the collisions with the external boundary ($M_{\out}$) and with the obstacle ($M_{\inn}$)}
\label{fig:dinamicas}
\end{figure}

It is important to have in mind that the circular billiard alone, without the inner obstacle, is completely integrable. 
On the other hand, the dynamics purely generated by the inner circle,  is somehow equivalent to Sinai's billiard 
which is ergodic.
Billiards in annular tables  
may share the properties of these two classical examples and exhibit a combined mixed dynamics. 
In particular, it was observed that the annular billiard undergoes very interesting dynamical bifurcations as one varies the parameters. 
The dynamics ranges from integrability (when the circles are concentric) to {\em chaotic} (when the obstacle is small and close to the exterior boundary).
Between these two extreme situations, the typical mixed Hamiltonian dynamics appears with elliptical islands surrounded by chaotic regions. The complexity of this dynamics, as observed numerically in \cite{saito}, can be seen on Figure~\ref{fig:dinamicas}. 

We have here a scenario similar to other families of conservative systems as, for instance,  the 
standard map \cite{sm,gorodetski}.
Among billiards we mention, on one hand, 
the different types of stadiums \cite{cmss,squash} and mushrooms \cite{mush} and,
on the other hand, the moon \cite{moon} and lemon \cite{lemon} billiards. 

More general systems of convex boundaries with inner scatters have been studied by several authors, among them we cite  Foltin,  Chen and Bolotin \cite{foltin,chen,bolotin}. 
Foltin \cite{foltin}  showed that, for a generic choice  of convex external boundary, the system has positive topological entropy provided the obstacle's radius is small enough.
This result was also obtained by Chen  in \cite{chen} using different techniques. 
In both works the dynamics around a specific class of periodic trajectories colliding orthogonally with the obstacle (we will refer to these as {\em normal periodic trajectories}) is at the center of the proof.
The result follows from the fact that under certain generic conditions the dynamics in an neighborhood of  such normal periodic orbits is conjugated to a shift. 
Bolotin \cite{bolotin} has proved that also under generic conditions on the external boundary and for small obstacles, the system presents hyperbolic sets around the normal periodic orbits. 
All these results were obtained by perturbation of the convex boundary and apply to generic situations. 
They do not apply to the annular billiard, as the the circular shape of the exterior boundary is certainly not generic. 
Nevertheless we use here some similar techniques and the role of normal trajectories is also central.  

Billiards in annular tables, as they have both convex and concave components, are examples of the so-called focusing and dispersing dynamics. 
It is well known that the dynamical behavior of such systems depend on the balance of these two effects. 
In systems with convex components sufficiently far apart
hyperbolicity is generated through the defocussing mechanism \cite{buni,donnay}. 
This is not the situation of the annular billiard and the standard defocussing arguments do not apply, as well as for moon and lemon billiards \cite{lemon,moon}. 
These three models fall in a different category of systems, where hyperbolicity comes from other mechanisms.
In this work, we show how to calibrate the distance between the centers and the radius of the inner circle in order to obtain hyperbolicity in large parts of the phase space of the annular billiard.
As pointed by previous results and numerical experiments, this hyperbolicity occurs for large eccentricity and small obstacle, i.e. $\delta \approx 1$ and $r \approx 0$.

Thus, if on one hand hyperbolicity seems to come from the dispersive obstacle, on the other hand the convexity of the outside boundary seems to be related to stability and KAM phenomena. 
The existence of stable periodic orbits in the annular billiard, in particular the one of period two, was investigated by  Saito et al. \cite{saito} 
{\ver
while an extensive study of other periodic orbits was presented by Gouesbet et al \cite{gouesbet}. 
For small eccentricity stable orbits of small period can be observed  \cite{rbrazm, malu} .
The large island of a period two trajectory which exists for $r>\delta$ is clearly visible on Figure \ref{fig:dinamicas} as well as a  period four and a period six. As one changes the parameters in the opposite direction by increasing the eccentricity (and  decreasing the size of the obstacle)  the islands become smaller and the orbits undergo transitions from elliptic to hyperbolic. 
In particular Saito et al. point out that the system seems to become ergodic as $\delta \to 1$ and $r \to 0$ if one considers only the trajectories which collide with the obstacle.  In this work we focus on this last situation and look for hyperbolic behavior.}  
More recently, Dettmann and Fain \cite{dettmann} have exhibited families of stable normal periodic orbits in the annular billiard when the obstacle is small and near the boundary, concluding that the system can not be ergodic for open sets of values of parameters close to this limit. 
The result is obtained through an explicit construction of suitable orbits and a direct computation of their non linear stability.
The existence of elliptical island follows from Birkhoff's Normal Form and Moser's Twist Theorem.
This is sometimes a tricky problem involving hard computations in very specific situations. 
Also for small obstacles and large eccentricity, we obtain elliptical islands associated to a  bifurcation of homoclinic tangencies
(Newhouse phenomenon).
As far as we know, this is the first time that this mechanism is described explicitly in billiards.

The presence of these elliptical islands gives a negative answer to the question if, as in the case of stadium like systems, there is a region in the parameter space where the annular billiard is fully chaotic (i.e. has positive Lyapunov exponent in a region of full measure).
It is natural to ask 
if there are any values of the parameters 
such that the chaotic region has positive measure. 
{\ver 
This is a challenging question, as well as the question of the size of the region occupied by islands. 
Another challenging question is the existence (or not) of other dynamical elements characteristic of conservative systems  such that invariant rotational curves and Aubry-Mather sets.}

The goal here is to present a global picture of the dynamics on annular billiards for large eccentricity and small obstacle
considering hyperbolic and non-hyperbolic properties. 
The big picture we obtain is  the following: 
for many values of parameters corresponding to a small eccentric obstacle, the system presents an ``almost dense" hyperbolic horseshoe, corroborating the numerical observation of ``chaos" in \cite{saito}. 
However, the  constructed hyperbolic set  has zero measure 
and, in many cases,  coexists with an also``almost dense" set of elliptical islands
originated from the generic bifurcation of quadratic homoclinic tangencies (Newhouse phenomenon).

Obtaining results for the annular billiard is somehow simplified by the fact that it is generated by two simple dynamics where the calculations can be made explicitly. What makes the situation more delicate is that the system is singular due to the existence of trajectories that are tangent to the obstacle which implies a loss of regularity and of compacity.
To overcome this difficulty, in Section \ref{sec:preliminares}, after writing down the billiard map explicitly, we describe the domain of the first return to the obstacle map $G_{\delta,r}$ with a special attention to the image and preimages of its boundary. Understanding the geometry of these singular curves as the parameters 
 $(\delta,r)$ converge to  $(1,0)$ is crucial along this work.

The key point to obtain hyperbolicity is that,
through a careful analysis of the tangent map $DG_{\delta,r}$ as $r \to 0$,
we can identify a strong expansion direction  in a certain region of the phase space. 
This fact gives rise to the hyperbolicity since it allows a construction of a cone field as described in Section~\ref{sec:cones}.
These cones are  preserved along orbits staying in these regions, which motivates the search of periodic orbits.  
In particular, a family of normal  periodic points is at the core of construction of a hyperbolic set and 
in Section~\ref{sec:normal} we show  the abundance of these points.

Our first result is the existence of hyperbolic sets  which become ``large", in the sense that they converge to the entire phase space, as the obstacle decreases in size and approaches the external boundary:
\begin{mtheorem} There is an open set of parameters $\Omega_0$ accumulating $(1,0)$
and a piecewise continuous family  $\Omega_0 \ni (\delta,r) \mapsto \Lambda_{\delta,r}$ of horseshoes for the first return to the obstacle map $G_{\delta,r}$ such that the maximum distance of  any point of the phase space  to $\Lambda_{\delta,r}$ goes to zero as $(\delta,r) \rightarrow (1,0)$. 
\label{thmA}
\end{mtheorem}

The proof of the above theorem is in Section \ref{sec:hyperbolic}. 
A compact invariant set in the region of hyperbolicity  is constructed from
the normal trajectories
which,  colliding orthogonally with the obstacle,
originate periodic points of the first return map.
{\ver 
Using the cones described in Section~\ref{sec:cones}  in  subsets around these periodic orbits
we construct hyperbolic invariant sets, the horseshoes mentioned in the statement of the theorem. 
As a part of this construction,
we obtain a symbolic description of the dynamics in the hyperbolic set $\Lambda_{\delta,r}$.
We show that the map $G_{\delta,r}$ restricted to it is conjugated to a subshift with a number of symbols which grows to infinity as  $(\delta,r) \to (1,0)$.
}

Analytically, a normal orbit corresponds to the intersection of two curves in the phase space,
which is transverse in the hyperbolic region.
Outside the hyperbolic region there are 
tangent normal points which turn out to be closely related to non-hyperbolicity.  
In the last section we show how tangent normal periodic points give rise to tangencies of invariant manifolds. 
We are able to show that for many parameters, quadratic homoclinic tangencies between manifolds of points in the set $\Lambda_{\delta,r}$ appear:

\begin{mtheorem} There is a set $\Omega'_0 \subset \Omega_0$ accumulating $(1,0)$ such that the maps  $G_{\delta,r}$, for $(\delta,r) \in \Omega_0'$ present  quadratic homoclinic tangencies unfolding generically with the parameter $r$.
\label{thmB}
\end{mtheorem}

Unlike the general setting of quadratic tangencies in dimension two between invariant  one-dimensional foliations, where one has to deal with the delicate analysis of intersection of Cantor sets \cite{cs}, here the reversibility of the system plays a major role.
This follows from the fact that if a branch of a stable manifold intersects the symmetry curve then we automatically obtain a {homoclinic} point of the basic set. 
Thus, quadratic tangency between  a branch of a stable manifold with the symmetry curve implies quadratic 
{homoclinic} tangency.

 As a consequence of the bifurcation of the homoclinic tangencies, the annular billiard presents the so-called Conservative Newhouse Phenomenon
with the appearance of {many} elliptical islands.
In fact, a detailed analysis of the bifurcation process enables one to use Duarte's theorem \cite{duarte}  in order to prove 
that, for many values of the parameters, the annular billiard has elliptical islands scattered across the phase space. 
This is the content of our third theorem, also in Section \ref{sec:tangencias}.
The last statement of the theorem strongly relies on our accurate description of the hyperbolic sets.

\begin{mtheorem} 
\label{thmC}
There is a set $\Omega_0''$ accumulating $(1,0)$ such that if $(\delta,r) \in \Omega_0''$ then  the map $G_{\delta,r}$ has a set $\mathcal{E}_{\delta,r}$ of generic elliptic periodic points. Moreover the distance of any point of the phase space to $\mathcal{E}_{\delta,r}$ tends to zero as $(\delta,r) \rightarrow (1,0)$.
\end{mtheorem}

In short, we were able to discriminate dynamical structures that appear in the phase space of the annular billiard (hyperbolic sets and elliptical islands) in the small and eccentric obstacle limit. 
However, we do not have an estimative of the measure of the 
{\ver chaotic region}
and in fact we don't even know if it is positive. 
Moreover, even if the
elliptical islands clearly sum up to a positive measure region, we do not know its extension.
{\ver
The estimates of the size of a specific elliptical island are in general hard to produce, as they usually involve a thorough analysis  of normal forms. Moreover, concerning the islands resulting from the Newhouse phenomenon in our case, we only know, as a general fact, that they exist and have long period (and so small islands).}

As a conclusion, we mention that 
there are several interesting questions concerning annular
billiards besides the (Lebesgue) 
measure of the chaotic region.
{\ver For instance, it is natural to ask if the closure of the union of the hyperbolic set(s) we produce and the islands has full measure and, if not, what is its complement?}
 Furthermore, one would like to have a more precise description of the bifurcation set, specially the parameter set corresponding to homoclinic or heteroclinic tangencies.
There is also the question of the dynamics inside an elliptical island, from the point of view of Zehnder's genericity or the existence of instability regions (in the sense of Birkhoff) containing an hyperbolic set inside an island.
This problem is possibly related to the destruction of invariant curves for parameters near the concentric case 
{\ver (although this is a very degenerate situation)} or to the transition of the stability of the orbits of period two (trajectories orthogonal both to the obstacle and to the external boundary)  or higher.
Finally, we point out that 
some of the results we obtained here are also true for generic external convex boundaries \cite{tese}.

To summarize, the sketch of the paper is the following:
In Section\ref{sec:preliminares}, we present the annular billiard and the first return to the obstacle maps as well as the domain, with a special attention to the singularities.
Section~\ref{sec:cones} contains the definition and properties of a cone field and Section~\ref{sec:normal} the description of normal periodic orbits.
These are the ingredients to construct a hyperbolic set and prove Theorem \ref{thmA} in Section~\ref{sec:hyperbolic}.
Finally, in Section~\ref{sec:tangencias} we show how  homoclinic tangencies are produced 
(Theorem \ref{thmB}) and, as a consequence, we have the existence of elliptical islands (Theorem \ref{thmC}).

%% file: 2-preliminares.tex
\section{Preliminaries}
\label{sec:preliminares}

{\ver
The billiard problem originally consists in the description of the free motion of a point particle in a bounded region of the plane with elastic collisions at the boundary. Conservation of energy and linear momentum implies the reflexion law at impacts. 
As a conservative system with two degrees of freedom,  each state is given by a point in the region and an unitary vector which accounts for the direction of motion. After some identifications, the time evolution is given by a three dimensional flow \cite{chernmark,foltin}, which in our case is defined for all time. It is usual to study the billiard dynamics through a restriction to the Poincaré section taken at the boundary of the region. The billiard map is then defined by the first return to boundary and thus associates to each impact, the next one.}

Given an annular region $Q_{\delta,r} \subset {\mathbb R}^2$ we assume that the normal vectors  point inside it. The external circular boundary $\gamma$  is parametrized by its central angle $s \in \mathbb{S}^{1}$ and is oriented counterclockwise, while the inner circular obstacle $\alpha$ is parametrized by its central angle $\omega \in \mathbb{S}^1$, and is oriented clockwise. 
Here we consider 
${\mathbb S^1} \sim (-\pi,\pi]$.
As usual, the billiard map is described by two variables: one for the position on the boundary ($s$ or $\omega$) and one for the direction of the trajectory, given the oriented angle from the inward normal vector to  the outgoing velocity ($\theta$ at the exterior boundary and $\beta$ at the obstacle). 
A collision with the external circle $\gamma$ is  then represented by a point $(s,\theta)$ in  the {open  cylinder} 
$M_{\out}= \mathbb{S}^1\times (-\frac{\pi}{2},\frac{\pi}{2})$  and a collision
with the obstacle $\alpha$  is represented by a point $(\omega , \beta)$ in the {closed cylinder} 
$M_{\inn}=\mathbb{S}^1\times [-\frac{\pi}{2},\frac{\pi}{2}]$.
The disconnected 
phase space of the billiard map  $T=T_{\delta,r}$ is the union $M=M_{\out} \cup M_{\inn}$.
{\ver We observe that the map may be extended by considering the boundary of $M_{\out}$ as fixed points.} 
In order to lighten our notation we will frequently omit the subscript $\delta,r$ that indicates the dependence of the maps and sets on the parameters. We will also refer to the inner circle as {\em the obstacle} and to the external circle simply
 as {\em the circle} or {\em the boundary}.

So, $T: M \rightarrow M$  denotes the billiard map in the annular region in general.
It is well defined and  is invertible, as it is reversible with respect to the involution $R(a,b)= (a,-b)$, i.e $ T^{-1}=R  \circ T \circ R$.
This reversibility implies that the phase space is symmetric with respect to the middle horizontal line. In particular, every orbit has its symmetrical which corresponds to the same trajectory traveled in the opposite direction.   {\em Symmetric orbits} are invariant and reverse orientation.

As described below, $T$ is  defined by parts: it is a piecewise diffeomorphism with a singular set generated by the tangent collisions with the concave obstacle.
$T$ is globally $C^0$ and piecewise $C^{\infty}$.

\begin{figure}[h]
\includegraphics[width=0.22\textwidth]{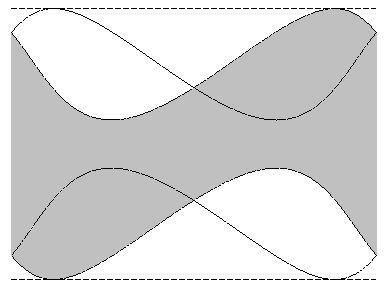}\hfill
\includegraphics[width=0.22\textwidth]{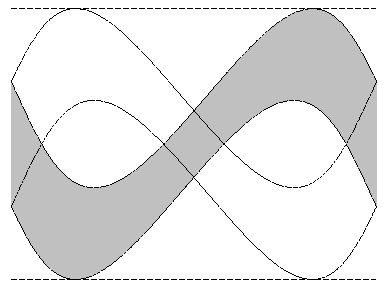}\hfill
\includegraphics[width=0.22\textwidth]{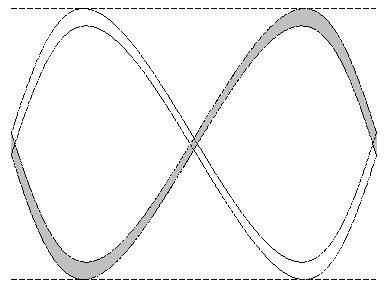}\hfill
\includegraphics[width=0.22\textwidth]{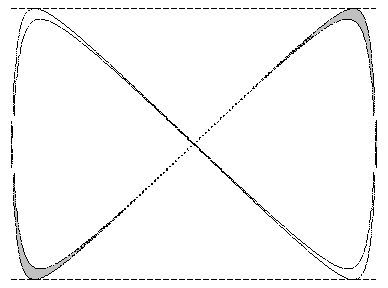}
\caption{$M_{\inn}^+$(gray) and $M^-_{\inn}$(white) in the complement of the whispering gallery in $M_{\out}$ (scaled), for $r>\delta$, $r<\delta$, $r \ll \delta$,
$r \ll \delta \approx 1$ }
\label{fig:malpha}
\end{figure}

To describe  $T$, we must distinguish between three different situations: the collisions from the obstacle to the (external) circle, from the circle to the obstacle and from the circle to the circle. 
Any trajectory from the obstacle will hit the circle in the sequence, 
{\ver which implies that $T(M_{\inn}) \subset M_{\out}$ and $T^{-1}(M_{\inn}) \subset M_{\out}$}.
 A trajectory leaving the circle will hit the obstacle if and only if 
$|\sin \theta + \delta \sin (\theta-s)| \le r$. We introduce the sets 
\begin{eqnarray} 
&T(M_{\inn}) &=  M_{\inn}^+=\{(s,\theta): |\sin \theta + \delta \sin (\theta+ s)| \leq r\} 
\label{eqn:M+-}
\\
&T^{-1}(M_{\inn}) &=  M_{\inn}^-= \{(s,\theta): |\sin \theta + \delta \sin (\theta- s)| \leq r\} \nonumber
\end{eqnarray}
which are topological cylinders in $M_{\out}$ ({\ver Figure~\ref{fig:malpha}}).

The restriction $T: M_{\inn}\rightarrow M_{\inn}^+$ (from the obstacle to the boundary) is implicitly given by 
\begin{eqnarray} T(\omega,\beta)=(s,\theta) \hbox{ with } \left\{\begin{array}{ll}
\sin \theta +\delta \sin(\theta+s)=-r\sin\beta&\\
          \omega+\beta= -s-\theta & \\
       \hbox{and }   |\sin\theta+\delta\sin(\theta+s)|\leq r&
\end{array}\right.
\label{eqn:T+}
\end{eqnarray}

Considering the trajectories leaving the exterior boundary,  the restriction $T: M_{\inn}^-\rightarrow M_{\inn}$ (from the boundary to the obstacle) is implicitly given by 
\begin{eqnarray} T(s,\theta)=(\omega,\beta) \hbox{ with } \left\{\begin{array}{ll}
\sin \theta +\delta \sin(\theta-s)=-r\sin\beta&\\
          \omega-\beta= \theta-s & \\
       \hbox{where }   |\sin\theta+\delta\sin(\theta-s)|\leq r&
\end{array}\right.
\label{eqn:T-}
\end{eqnarray}

In the particular case of a trajectory from the boundary to the boundary without  colliding with the obstacle the map $T:M_{\out}\setminus M_{\inn}^- \rightarrow M_{\out}$ is given by the circular billiard map (denoted by $F$)
\begin{eqnarray}
T(s,\theta) = F(s,\theta)=(s + \pi -2\theta, \theta) \label{eqn:F}
\end{eqnarray}

The map from the exterior boundary is then clearly discontinuous on $T^{-1}(\partial M_{\inn})$. 
The concavity of the obstacle implies so the existence of unavoidable  tangent collisions corresponding to $\partial M_{\inn} = \{ \beta = \pm \pi/2 \}$.  Moreover, the map $T$ is not differentiable on $\partial M_{\inn}$.
The map $T$ has then a singular set given by the curves $\partial M_{\inn} \cup T^{-1}(\partial M_{\inn})$ and the inverse $T^{-1}$ has a singular set $\partial M_{\inn} \cup T(\partial M_{\inn})$. 
Out of the singular set, $T$ is a $C^\infty$ diffeomorphism.

It is also well known that the billiard map is conservative and preserves the measure $\mu$ given in $M_{\out}$ by 
$d\mu = \cos \theta ds d\theta$ and in $M_{\inn}$ by $ d \mu = r\cos\beta d\omega d\beta$.
{\ver In the canonical variables (tangential momentum and arc length) the Lebesgue measure $dpds$ is preserved.} 
This can also be directly checked from the expressions \ref{eqn:T+}, \ref{eqn:T-} and \ref{eqn:F} above.
The choice of the coordinate $\omega$ in $M_{\inn}$ instead of the usual arc length is particularly convenient as we want to use arguments with $r \to 0$.

 
Trajectories leaving the exterior boundary in almost tangential directions will circulate around without hitting the obstacle. More precisely, if a trajectory leaves the exterior boundary with an angle  $|\sin \theta| > \delta +  r$, it will follow the circular billiard motion forever 
{\ver with a circular caustic concentric to the boundary. The corresponding invariant region in $M_{\out}$ is a cylinder foliated by invariant rotational horizontal curves, clearly visible on Figure~\ref{fig:dinamicas}.
This region is known as the {\em whispering gallery}  and we denote it by $M_{w} \subset M_{\out}$.
Besides the whispering gallery there are other trajectories from the exterior boundary which do not hit the obstacle. They necessarily correspond to periodic orbits in
$M^c_w \backslash (M^+_{\inn} \cup M^-_{\inn}) \subset M_{\out}$ as any non periodic trajectory in the circular billiard would be dense on the caustic 
of radius $|\sin \theta| < r+ \delta $ and so cannot avoid the obstacle.  
}

On the other hand, a trajectory leaving the obstacle and hitting the boundary with {$|\sin \theta| < \delta +  r$} 
will hit the obstacle again an infinite number of times. 
{\ver
As for $|\sin \theta| = \delta +r$ 
the only trajectories leaving the obstacle and not returning to it
occur when the orbit of a tangent point $(\omega,\beta) = (\pi,  \pm \frac{\pi}{2})$ of $M_{\inn}$ is not periodic
(the trajectory corresponds to the circular billiard caustic of radius $\delta + r$).
So, with a possible exception of two points, every point }
$(\omega,\beta) \in  M_{\inn}$ has a finite return time to $M_{\inn}$ 
$$\nu (\omega,\beta)= min\{j \geq 1: T^j(z) \in M_{\inn}\}$$
This will allow us to define  the first return to the obstacle map
 \begin{equation*}G = M_{\inn} \rightarrow M_{\inn}, \,\,\,\, 
G(\omega,\beta)= T^{\nu(\omega,\beta)} (\omega,\beta)= T\circ F^{\nu(\omega,\beta) - 2} \circ T (\omega,\beta)\end{equation*}
So, a trajectory leaving the obstacle from $(\omega_0,\beta_0)$ will return to  it at 
 $(\omega_1,\beta_1) = G  (\omega_0,\beta_0) $  after $m = \nu(\omega_0,\beta_0)-2$ collisions with the external boundary at points $(s_0,\theta_0) = T(\omega_0,\beta_0),  \ldots ,  (s_m,\theta_0) = F^m (s_0,\theta_0)$ with $(\omega_1,\beta_1) = T (s_m,\theta_0)$.

From the properties of the billiard map, $G$ is a piecewise $C^\infty$ diffeomorphism. 
We denote its set of singularities  by $\mathcal{S}^- = \partial M_{\inn} \cup G^{-1}(\partial M_{\inn})$. 
The singular set of $G^{-1}$ is denoted by  $\mathcal{S}^+ = \partial M_{\inn} \cup G(\partial M_{\inn})$. 
 


The annular billiard has two period 2 trajectories, bouncing between the obstacle and the exterior boundary with orthogonal collisions. They correspond to the fixed points of the first return map $G$: $(0,0)$ and $(\pi, 0)$. The second is always hyperbolic, while the first one is hyperbolic if $r<\delta$ and elliptic (in fact, Moser stable) if $r>\delta$
\cite{saito,malu,rbrazm,rbraz}. 
The annular billiard has also many other periodic normal trajectories, as  the orbits presented in Section \ref{sec:normal}, which play a very important role in the dynamics and in our analysis, as in
\cite{foltin,chen}. The stability of some of these orbits was established in \cite{dettmann}. 
It is clear that the stability of periodic orbits and so the dynamics depend on the parameters. 
{\ver In particular, numerical experiments seem to indicate that, besides the period two orbit, the other short period elliptical orbits also loose stability as $r$ decreases. }
This is one reason why, in order to investigate the chaotic behavior of the annular billiard, 
we focus in the dynamics for small $r<\delta$.
More precisely, we will present results on the two parameter family of maps $G_{\delta,r}$, describing some aspects of the dynamics as $(\delta,r) \rightarrow (1,0)$.  

In our strategy, the parameter dependence of some relevant subsets of the phase space is very important.
The sets $M^+_{\inn}$ and 
$M^-_{\inn}$ are contained in the cylinder $ M^c_w = |\sin \theta | \le r + \delta$, the complement of the whispering gallery in $M_{\out}$.
The boundaries of these cylindrical sets
as defined in \ref{eqn:M+-}
are given by the curves
$$ \partial M^{\pm}_{\inn} = |\sin \theta + \delta \sin (\theta \pm s) | = r $$ 
{\ver which have a single point of tangency with the top and the bottom of $M^c_w$.
For fixed $\delta$,  the whispering gallery grows when the obstacle becomes smaller and so these sets  becomes thinner as $r$ goes to 0, as we will precise bellow (see Figure~\ref{fig:malpha}).}

If we denote the horizontal line corresponding to orbits leaving the obstacle in the normal direction by
\begin{equation}
L^0  = \left\{(\omega, \beta) \, | \, \beta = 0  \right\}\subset M_{\inn}
\label{eqn:L0}
\end{equation} 
its image and preimage in $M_{\out}$ are defined by
\begin{eqnarray}
 L_{\delta}^+ =&T(L^0) &= \{(s,\theta): \sin \theta + \delta \sin (\theta+ s)= 0\}  
\\
L_{\delta}^-=& T^{-1}(L^0) &=\{(s,\theta): \sin \theta + \delta \sin (\theta- s)=0\} \nonumber
\end{eqnarray}
We notice that these sets depend only on $\delta$, the eccentricity parameter.

It follows that the boundaries $\partial M_{\inn}^+$ and  $\partial M_{\inn}^-$  converge as $r$ goes to 0 respectively to the curves $ L_{\delta}^+$ and $L^-_{\delta}$.
Therefore, as $r \to 0$,  the subsets $M_{\inn}^\pm$ 
become narrow cylindrical strips also converging to the curve
$L_{\delta}^\pm$. This contraction has deep consequences on the dynamical behavior.
We also point out that  the curves  $ L_{\delta}^+$ and $L^-_{\delta}$ are graphs of analytic functions 
of $\theta$ converging uniformly in $(-\pi, \pi)$, as $\delta \to 1$ to the lines $2 \theta \pm s = 0$. 
{\ver As for any $\delta$, $s=\pm \pi$ implies $\theta = 0$, the limit is strongly discontinuous.}

Another relevant  preliminary observation  is that  for  $\delta > r$,
the domains $M^+_{\inn}$ and  in $M^-_{\inn}$ do not contain
any horizontal line $\theta = constant$. Moreover, in this case, the intersection $M^+_{\inn} \cap M^-_{\inn}$ as two distinct connected components, each one  containing one period two orbit, corresponding to the two fixed points  of the first return to the obstacle map: $(0,0)$ and
$(\pi,0) \in M_{\inn}$. 

{\ver We note again that,  in order to make the notation lighter and the reading easier, we will drop the subscripts $\delta,r$ of maps and sets in our proofs and computations whenever the parameters are fixed and the dependence on them is clear.}

%% file: 3-cones.tex
\section{Finding regions of hyperbolicity: Cone Fields}
\label{sec:cones}
In this section we will show that the annular billiard presents hyperbolicity
for a {wide} choice of parameters.
{
This hyperbolicity follows from the existence of a cone field, which, in some region of the collision with the obstacle set $M_{\inn}$, is strictly preserved by the first return to the obstacle map $G$. 
As vectors in the cone are uniformly expanded, any invariant compact set will be uniformly hyperbolic. 
}

Defining horizontal/vertical cone fields $z \mapsto C^{\pm}(z)$ for  $z=(\omega,\beta) \in int \, M_{\inn}$ 
by
\begin{eqnarray} 
C^+(z) := \{u=(u_1,u_2) \in {\cal T}_{z}M: u_2.u_1 \geq 0\} 
\\
C^-(z) :=  \{u=(u_1,u_2) \in {\cal T}_{z}M: u_2.u_1 \leq 0\} \nonumber
\label{def:cones+-}
\end{eqnarray}
we have
\begin{theorem} \label{Thm:Cones+-}
There is a subset $\Omega_* \subset \Omega$ of parameters,
such that for each $(r,\delta) \in \Omega_{*}$ there are subsets $H_{\delta,r}^{\pm}  \subset M_{\inn}$
with  $H_{\delta,r}^+=G_{\delta,r}(H_{\delta,r}^-)$ where
\begin{enumerate}
\item[(i)] 
the map $ \displaystyle G_{\delta,r}: H_{\delta,r}^- \rightarrow H_{\delta,r}^+$ (resp. $G_{\delta,r}^{-1}: H_{\delta,r}^+ \rightarrow H_{\delta,r}^-$) strictly preserves  the cone field $C^+$(resp. $C^-$).
\item[(ii)]  for points in $H_{\delta,r}^{-}$ (resp $H_{\delta,r}^{+} $)
\begin{equation} 
||D G_{\delta,r} ||  \hbox{ (resp. $||D G_{\delta,r}^{-1}||$ )}  \geq \rho  \label{expand}
\end{equation}
with $\rho \rightarrow \infty$ as $r \rightarrow 0$.
\end{enumerate}
\end{theorem}

As a consequence
\begin{corollary} If $(\delta,r) \in \Omega_*$ and $\Lambda \subset H_{\delta,r}^-$ is a compact invariant set for $G_{\delta,r}$
 then $\Lambda$ is a uniformly hyperbolic set for $G_{\delta,r}$.
\label{Coro:Cones}
\end{corollary}

As usual, we will drop the subscript $\delta,r$ in maps and sets.

Let us consider a trajectory leaving the obstacle with $(\omega_0,\beta_0) \in M_{\inn} \backslash {\mathcal S}^-$
and returning to it with $(\omega_1,\beta_1)=G(\omega_0,\beta_0)$, after $m+1$ impacts with the exterior border  $\gamma$ given by $\{(s_0,\theta),....,(s_{m},\theta)\}$. 
{A straightforward computation from equations  \ref{eqn:T-}, \ref{eqn:T+} and \ref{eqn:F} leads to following expression of the derivative of the map}
\begin{equation}
DG{(\omega_0,\beta_0)}= \left(
\begin{array}{cc}
         a_{11} & a_{12}\\
          a_{21}&a_{22}  \end{array}\right)= a_{21} \left(
\begin{array}{cc}
         1 & 1\\
          1&1  \end{array}\right) + 
					\left( \begin{array}{cc}
         \widetilde{a}_{11} & \widetilde{a}_{12}\\
          0 & \widetilde{a}_{22} \end{array}\right)
          \label{eqn:DT}
\end{equation}        
where 
\begin{equation}
 a_{21} = 
-\frac{\cos \theta}{r \cos \beta_1}
\left(
\frac{\delta \cos \varphi_0}{\cos \theta} + \frac{\delta \cos \varphi_1}{\cos \theta} 
+ 2(m+1) \frac{\delta \cos \varphi_0}{\cos \theta} \, \frac{\delta \cos \varphi_1}{\cos \theta} 
\right) 
\label{eqn:a21}
\end{equation}
and
\begin{eqnarray}
 \widetilde{a}_{11} & = &  \displaystyle 1 + {2(m+1)}\frac{\delta\cos\varphi_0}{\cos\theta} \label{eqn:atilde} \\
\widetilde{a}_{22} & = &  \displaystyle \frac{\cos\beta_0}{\cos\beta_1} \left(  1 + {2(m+1)}\frac{\delta\cos\varphi_1}{\cos\theta}\right) \nonumber\\
 \widetilde{a}_{12} & = & 
\displaystyle 1+ \frac{\cos\beta_0}{\cos\beta_1} 
+ 2(m+1) \left(\frac{\delta\cos\varphi_0}{\cos\theta} + \frac{\cos\beta_0}{\cos\beta_1}\frac{\delta \cos\varphi_1}{\cos\theta}- r\frac{\cos\beta_0}{\cos\theta} \right)\nonumber\\
& = & \widetilde{a}_{11} + \widetilde{a}_{22} - 2 r (m+1) \frac{\cos\beta_0}{\cos\beta_1}
\nonumber 
\end{eqnarray}
Here $\varphi_0 = s_0 + \theta = -\omega_0 -\beta_0$ and $\varphi_1= s_m - \theta = -\omega_1 + \beta_1$ represent  the angle between the outgoing trajectory leaving the obstacle (respectively incoming  back) and the horizontal direction.
It is worthwhile to note that $\det DG = \cos \beta_0 /\cos\beta_1$.

The key observation is that as $r$ approaches zero, the first  matrix in the sum dictates the behavior of the tangent map as long as $a_{21} \ne 0$.
However, it is easy to check that $a_{21}$ is negative at $(0,0)$  and  positive
at $(\pi,0)$, corresponding to the 2-periodic trajectories. On the other hand,  it is clear that $a_{21}$ vanishes
 if the trajectory paths between the obstacle and the boundary are vertical,  as $\varphi_0 = \varphi_1 =  \pm \pi/2$. These observations indicate that there is no hope to bound $a_{21}$ away from zero globally on $M_{\inn}$. 
Our strategy is then  is to find  a subset of parameters {$\Omega_*$} 
and a subset of phase space ${H}^-_{\delta,r}$ 
where all the entries of the matrix $DG$ are non zero and have the same sign. This will imply  the preservation of the cone $C^+$ for $G$ and, by reversibility, also implies the preservation of the cone $C^-$ by $G^{-1}$ \cite{woj1,woj2}.

\begin{lemma} Let $\displaystyle \zeta = \frac{\delta  \cos \varphi}{\cos \theta}
$. If $\delta^2 > \frac{1}{2} $ and $ r< \frac{1}{4}(\delta-\delta^2)$ then for any $\varphi \in [0,2\pi]$, $\theta \in [-\pi/2,\pi/2]$ such that  $|\sin \theta + \delta \sin \varphi | \le r $ and $| \sin \theta | \le \delta^2$ we have
\begin{equation*}
\zeta_{min} =  \frac{\delta}{2}\sqrt{\frac{3}{1+\delta^2}}  < |\zeta| < \sqrt{\delta}= \zeta_{max} 
\end{equation*}
Moreover $\zeta_{min} >1/2$ and $\zeta_{max}  <  1$.
\label{lem:zeta0}
\end{lemma}
\begin{proof} If we use the coordinates $x = -\delta \sin\varphi$ and  $y = \sin \theta$, we have that 
$\zeta^2 = \frac{\delta^2 - x^2}{1-y^2}$ should be bounded on the compact parallelogram
$\{(x,y): |y-x|\leq r \,\, {\rm and} \, \, |y|\leq \delta^2\}$.
As $\nabla \zeta^2 = \frac{2}{1-y^2} (-x,\zeta^2 y)$, the origin is the only critical point inside the domain.
It is a saddle with $\zeta^2(0,0)=\delta^2$, and so minimum and maximum should be on the boundary.
Because of the symmetry of the function $\zeta^2$ we can restrict our search for the maximum and minimum values to the region bounded by the lines $y = x+ r$,
$y=x-r$, $y=\delta^2$ and the axis $x=0$, $y= 0$.  

The level curves of $\zeta^2 = k^2$ are the hyperbolas 
$
k^2 - \delta^2 = k^2 y^2 - x^2
$
and so, corresponding to $\zeta^2 = \delta^2$, we have the  asymptotes $x^2 = \delta^2 y^2$. 
The hyperbolas with vertices on the $y$ axis have $\zeta^2 > \delta^2$ and the ones with vertices on the $x$ axis correspond to $\zeta^2 < \delta^2$. 
This implies that the maximum value of $\zeta^2$ occurs on the segment of the line $y=x+r$ between the $y$-axis and the asymptote $x=\delta y$, i.e.,  between the points $(0,r)$ and $(\delta r/(1-\delta), r/(1-\delta))$ . 
Moreover, at the maximum point $(x^*,y^*)$ we have that $\nabla \zeta^2 . (1,1)$ = 0 and so the maximum value  
$\zeta^2(x^*,y^*) = x^*/y^*$. As $x^* = y^* - r$ if follows that   $\zeta^2 (x^*,y^*) = 1 - r/y^*$ and since $r<y^*<r/(1-\delta)$ we have that 
$$\zeta^2(x^*,y^*)< \delta < 1$$
Since the slope of the components of the boundary is 1 or 0, it is clear that no hyperbola with $\zeta^2 < \delta^2$ can have a tangency with them and so the minimum value must occur at a vertex. Comparing the values, and using that 
$r< 1/4  (\delta-\delta^2)$ is is easy to check that the minimum value is
$$
 \frac{\delta^2-(\delta^2+r)^2}{{1-\delta^4}} 
= \frac{((\delta-\delta^2) - r)(\delta+\delta^2+r)}{(1+\delta^2)(1-\delta^2)}> \frac{3}{4} \, \frac{\delta^2}{1+\delta^2} > \frac{1}{4}
$$
\end{proof}

\begin{figure}[h] 
\includegraphics[width=0.22\textwidth]{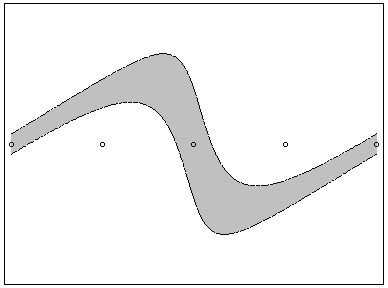}\hfill
\includegraphics[width=0.22\textwidth]{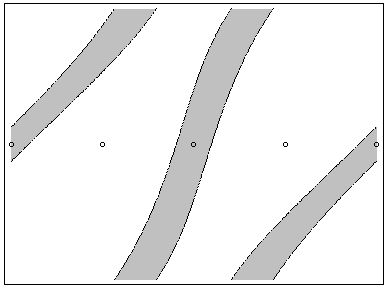}\hfill
\includegraphics[width=0.22\textwidth]{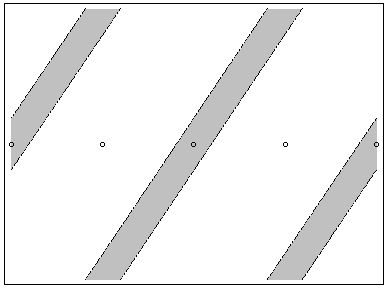}\hfill
\includegraphics[width=0.22\textwidth]{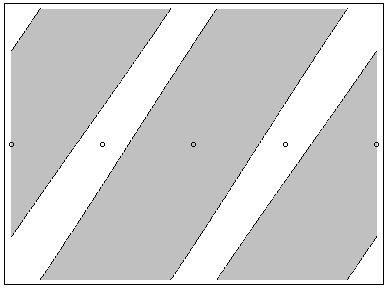}
 \caption{$H^+_{\delta,r} \subset M_{\inn}$ for $r>\delta$, $r<\delta$, $r \ll \delta$, $r \ll \delta \approx 1$  }
 \label{fig:H}
\end{figure}

Following Lemma \ref{lem:zeta0} above,  we define the horizontal strip 
$H_{\delta} = \{(s,\theta)  \hbox{ s.t. } |\sin \theta| < \delta^2 \} \subset M_{\out}$ and the subsets of $M_{\inn}$  
\begin{equation}
H_{\delta,r}^- =  T^{-1} (H_{\delta})  \ \ , \ \ H_{\delta,r}^+ = T(H_{\delta}) 
\hbox{ \ and \ } \bar{H}_{\delta,r}= H_{\delta,r}^+ \cap H_{\delta,r}^- 
\label{eqn:H}
\end{equation}
Its easy to check that $G(H_{\delta,r}^-)= H^+_{\delta,r}$ and $R(H^-_{\delta,r})= H^+_{\delta,r}$  and from equations $\ref{eqn:T-}$ and $\ref{eqn:T+}$ we have 
\begin{equation}
H^{\pm}_{\delta,r}= \{(\omega,\beta): |\delta\sin(\omega\mp \beta)+r\sin\beta|<\delta^2\}
\label{eqn:H+-}
\end{equation}
As noticed in the end of Section \ref{sec:preliminares}, if $r<\delta$,  the intersection of any horizontal strip in $M_{\out}$ with either $M^+_{\inn}$ or 
$M^-_{\inn}$ has two distinct connected components. So, 
 $H^\pm_{\delta,r} \subset M_{\inn}$ also have two connected components, one containing the point  $(0,0)$ and the other one the point  $(\pi,0)$. These components are bounded by the four curves with endpoints in $\partial M_{\inn}$ given by $ |\delta\sin(\omega\mp \beta)+r\sin\beta|=\delta^2$. Each one of the components of 
$M_{\inn} \backslash H^\pm_{\delta,r}$ contains one of the points $\left(-\frac{\pi}{2},0 \right)$ or $\left(\frac{\pi}{2},0 \right)$ (Figure \ref{fig:H}).

Fixing $\delta$ and taking $r\rightarrow 0$, the curves in $\partial H_{\delta,r}^- \backslash \partial M_{\inn}$ converge to the straight lines given by  $|\sin(\omega+\beta)|=\delta$.  Thus as $(\delta,r) \rightarrow (1,0)$, the components of $M_{\inn} \backslash H_{\delta,r}^-$ shrink to the  decreasing lines $\omega+\beta= \pm\frac{\pi}{2}$ and hence the set $H^-_{\delta,r}$  converges, in Hausdorff sense, to $M_{\inn}$. Similarly, as $(\delta,r) \rightarrow (1,0)$, the components of $M_{\inn} \backslash H^+_{\delta,r}$ shrink to the  lines $\omega-\beta= \pm\frac{\pi}{2}$ and $H^-_{\delta,r}$ converges to $M_{\inn}$.

\begin{lemma} For $\delta^2 > \frac{1}{2}$ and  $r< \frac{1}{4}(\delta-\delta^2)$, if $(\omega_0,\beta_0) \in H^-_{\delta,r}$ then
 $\left|a_{21}\right| \geq \frac{4A}{\sqrt{r}}$ where $A$  is a constant depending only on $\delta$.
\label{lem:a21}
\end{lemma}
\begin{proof}
Using the notation of Lemma \ref{lem:zeta0}, we can write
$$
a_{21} = - \frac{\cos \theta}{r \cos \beta_1} (\zeta_0 + \zeta_1+ 2(m+1) \zeta_0 \zeta_1) 
$$ 
As $(\omega_0,\beta_0)  \in \bar{H}$ and we have $r<\frac{1}{4}(\delta-\delta^2)$
$$\cos^2 \theta  > 1 - \delta^4 = (1+\delta^2)(1+\delta)(1-\delta) \ge \delta(1-\delta)  > 4r
$$
For any $m\geq 0$ we consider
$$g_{m} (x,y) =  2(m+1) x y + x + y$$
in the region $D=\{(x,y): \zeta_{min} \leq |x|, |y|  \leq  \zeta_{max} \}$, which corresponds to four equal squares in the plane and we  want to estimate its minimum value.
$g_m$ has a saddle point at $\left( \frac{-1}{2(m+1)} ,  \frac{-1}{2(m+1)} \right)$ with $g_m = -1/2(m+1)$. The level curves are hyperbolas with asymptotes through the saddle point parallel to the $x,y$ axis. There are two distint level curves with $g_m=0$, one through $(0,0)$ and the other one through $\left( \frac{-1}{m+1} ,  \frac{-1}{m+1} \right)$. These level curves are outside the four squares as $\zeta_{min} >1/2$ and $\zeta_{max} <1$ and so the minimum value should be on one of the corners. It is easy to check that in fact the minimum occurs at the vertex closest to the point  $\left( \frac{-1}{m+1} ,  \frac{-1}{m+1} \right)$. 
For $m=0$ it is the point $(-\zeta_{min},-\zeta_{min})$ while for $m\ge 1$ it is   $(-\zeta_{max},-\zeta_{max})$ .  
It follows that 
\begin{equation*}
| \zeta_0 + \zeta_1 + 2(m+1) \zeta_0 \zeta_1 |  \ge 
\left\{
\begin{array}{ll}
2 (\zeta_{max} - \zeta_{max}^2) & \hbox{ if $m=0$ }  \\
2 ((m+1) \zeta_{min}^2 -  \zeta_{min}) \ & \hbox{ if $m \ge 1$ }
\end{array}
\right .
\end{equation*}

So we have $|g_0| \ge  2(\sqrt{\delta} -  \delta)$ 
and for $m\ge 1$ 
$ \displaystyle |g_m| \ge |g_1|  \ge 2 \left( \frac{6\delta^2}{2(1+\delta^2)}- \frac{\delta\sqrt{3}}{2\sqrt{1+\delta^2}}\right)
$.
If we denote 
\begin{equation} 
A= \min \left\{(\sqrt{\delta}-\delta),  \left( \frac{6\delta^2}{2(1+\delta^2)}- \frac{\delta\sqrt{3}}{2\sqrt{1+\delta^2}}\right)\right\}
\label{eqn:constA}\end{equation}
the lower bound on $|a_{21}|$ follows.
In fact we have we have $g_0 \ge 2A$ and for $m \ge 0$
\begin{equation*}
|g_{m}| \geq (m+1) A
\end{equation*}
\end{proof}

\begin{lemma} For $\delta^2 > \frac{1}{2}$ and  $r< \frac{1}{4}(\delta-\delta^2)$, if $(\omega_0,\beta_0) \in H^-_{\delta,r}$  then 
\begin{equation*}
\left | \frac{\widetilde{a}_{11}}{a_{21}} \right | \le A_{11} \sqrt{r} \hskip 1cm
\left |  \frac{\widetilde{a}_{22}}{a_{21}} \right | \le A_{22} \sqrt{r} \hskip 1cm
\left | \frac{\widetilde{a}_{12}}{a_{21}} \right |\le A_{12} \sqrt{r} 
 \end{equation*}
where the constants $A_{11}, A_{22}$ and $A_{12}$ depend only on $\delta$.
\label{lem:aij}
\end{lemma}
\begin{proof}
Following the notation and definitions in the proof of Lemma \ref{lem:a21} above, we have
\begin{eqnarray*} 
\left|\frac{\widetilde{a}_{11}}{a_{21}}\right| & = & \left|\frac{r\cos\beta_1}{ \cos \theta}\right| |f_{m}(\zeta_0,\zeta_1)|\\  
\left|\frac{\widetilde{a}_{22}}{a_{21}}\right| & = & \left|\frac{r\cos\beta_0}{\cos \theta}\right| |f_{m}(\zeta_1,\zeta_0)| \\
\left|\frac{\widetilde{a}_{12}}{a_{21}}\right| &\leq  & \left|\frac{\widetilde{a}_{11}}{a_{21}}\right|+ \left|\frac{\widetilde{a}_{22}}{a_{21}}\right| +  \left| \frac{r\cos\beta_0}{\cos\theta}\right| \left| \frac{r\cos\beta_1}{\cos\theta}\right| \left| \frac{1}{g_{m}(\zeta_0,\zeta_1)}\right|
\end{eqnarray*}
where
\begin{equation*}  f_{m}(\zeta_0,\zeta_1)= \frac{1+2(m+1) \zeta_0}{g_{m}(\zeta_0,\zeta_1) }
\end{equation*}
If $ \zeta_{min} \le |\zeta_i| \le \zeta_{max}$
we have 
$$|f_m (\zeta_0,\zeta_1)|  \leq \frac{4(m+1)\sqrt{\delta}}{(m+1)A}= \frac{4\sqrt{\delta}}{A}$$
and  as $1/ \cos \theta < \sqrt{r}/2$ it follows that
\begin{equation*}
\left|\frac{\widetilde{a}_{11}}{a_{21}}\right| \leq \frac{2\sqrt{\delta}}{A} \sqrt{r}
\hspace{20pt} 
\left|\frac{\widetilde{a}_{22}}{a_{21}}\right| \leq \frac{2\sqrt{\delta}}{A}\sqrt{r}
\hspace{20pt} 
\left|\frac{\widetilde{a}_{12}}{a_{21}}\right| \leq \frac{4\sqrt{\delta}+1/4}{A}\sqrt{r}
\end{equation*}
\end{proof}

\begin{proof}[Proof of Theorem \ref{Thm:Cones+-}]
With the constant $A$ defined by equation~\ref{eqn:constA}, we consider the continuous function  
$$r(\delta) < \min \left \{ \frac{1}{4}(\delta-\delta^2), \frac{A^2}{(1/4+ 4\sqrt{\delta})^2}\right\}$$
and
define the set of parameters
\begin{equation}
\Omega_{*}= \{(\delta,r):\delta^2 > \frac{1}{2} \,\, {\rm and} \,\, 0<r<r(\delta)\}
\label{eqn:omega*}
\end{equation}
which is has no empty interior and accumulates $(1,0)$ as $A\to 0$ when $\delta \to 1$.

It is then clear  that if $(\delta,r) \in \Omega_{*}$ and $(\omega,\beta) \in H^-_{\delta,r}$, the matrix
$DG(\omega,\beta)$ as given in \ref{eqn:DT} has either positive or negative entries and  so the cone field $C^+$ is strictly preserved.
Moreover, taking $u=1/\sqrt{2}(1,1) \in C^+$ we have
$$
DG(\omega,\beta)  u =  a_{21} \left[
\left(
\begin{array}{c}
          \sqrt{2}\\
          \sqrt{2}
\end{array}
\right)
+
\frac{1}{\sqrt{2} \, a_{21}} \left(
\begin{array}{c}
          a_{11} + a_{12} \\
          a_{22}
\end{array}
\right)
\right]
$$
with, from Lemma \ref{lem:aij},  
$
e = 
\frac
{\sqrt{(a_{11}+a_{12})^2 + a_{22}^2}}
{ \sqrt{2} |a_{21}| }  < K \sqrt{r}
$ for some constant $K$ depending only on $\delta$.

We have then
\begin{equation}
|| DG(\omega,\beta)|| \ge || DG(\omega, \beta) \, u ||  \ge |a_{21}| ( 2- e)  > \rho 
\label{eqn:rho}
\end{equation} 
where $\rho$ is a constant depending only on $\delta$ which can be chosen using the bound on $a_{21}$ from Lemma \ref{lem:a21}. Moreover, $\rho \to \infty$ as $ r\to 0$.

By reversibility, $G^{-1}: H_{\delta,r}^+ \backslash {\mathcal S}^+ \rightarrow H_{\delta,r}^- \backslash {\mathcal S}^-$ preserves the cones $C^-$, expanding  vectors by the rate $\rho$.
\end{proof}

A closer look at the tangent map $DG$ as given by  \ref{eqn:DT} shows that as $r\to 0$ it strongly contracts vectors to the diagonal $(1,1)$ direction while the inverse $DG^{-1}$ contracts to $(-1,1)$.
We can use this fact to obtain more precise estimates on the expansivity and control on the hyperbolicity. To do so, we introduce the notion of stable and unstable curves, which play a fundamental role in our geometric arguments to exhibit both hyperbolic and non-hyperbolic behavior
in the annular billiard as studied in Sections $\ref{sec:hyperbolic}$ and $\ref{sec:tangencias}$.

\begin{definition}
For $(\delta,r) \in \Omega_*$, let
\begin{equation*}
c_1 =\min_{(\omega,\beta) \in {H}^{-}_{\delta,r} } \frac{a_{21}(\omega,\beta)}{a_{11}(\omega,\beta)} 
 \ \rm { and } \ 
c_2 = \max_{(\omega,\beta) \in {H}^{-}_{\delta,r}} \frac{a_{22}(\omega,\beta)}{a_{12}(\omega,\beta)}
\end{equation*}
and note that $c_1, c_2 \to 1$ as $r\to 0$.
A ${\mathcal C}^1$-curve   
$\ell(t) = (\omega(t),\beta(t))$ 
is called \textbf{unstable}  if   $\ell(t) \subset {H}^{-}_{\delta,r}$ and $c_1 \leq \frac{\beta'(t)}{\omega'(t)} \leq c_2$ 
and  it 
is called \textbf{stable} if $\ell(t) \subset {H}^+_{\delta,r}$ and $-c_2 \leq  \frac{\beta'(t)}{\omega'(t)} \leq -c_1$.
If  $\beta'(t)=0$  the curve is \textbf{horizontal}.  
\end{definition}


We summarize in the following propositions, some properties of stable and unstable curves which can be easily derived from the 
arguments leading to Theorem  \ref{Thm:Cones+-}

\begin{proposition} 
\label{prop:us-curves}
Let $\ell$ be a ${\mathcal C}^1$-curve and $|\ell|$ denote its length. 
Then, for  $(\delta,r) \in \Omega_*$ we have
\begin{enumerate}
\item
If $\ell$ is stable (unstable), then it is the graph $\omega=f(\beta)$ of a $1/c_1$-Lipschitz  monotone function $f$.
Moreover,  for $r \approx 0$ a stable (unstable) curve is $\mathcal{C}^1$ close to a straight segment of slope $-1$ ($1$).
\item If $\ell  \subset H^-_{\delta,r}$ is {either an unstable or a horizontal} curve then $G_{\delta,r}(\ell)$ is unstable and $|G_{\delta,r}(\ell)| \geq \rho |L|$.
\item If  $\ell \subset H^+_{\delta,r}$ is {either a stable or a horizontal} curve then $G_{\delta,r}^{-1}(\ell)$ is stable  and $|G_{\delta,r}^{-1}(\ell)| \geq \rho |\ell|$.
\end{enumerate}
As defined by \ref{eqn:rho},  the expansion rate $\rho$ depends only on $\delta$ and goes to $\infty$ as $r\to0$.
\end{proposition}

We also have the following description of the set of the singularities of $G: M_{\inn} \rightarrow M_{\inn}$, given by 
${\mathcal S}_{\inn}^-= \partial M_{\inn} \cup G^{-1}(\partial M_{\inn})$ and its inverse. 

\begin{proposition}  For $(\delta,r) \in \Omega_{*}$,
the singular set of the restriction $G_{\delta,r}: H_{\delta,r}^- \rightarrow H_{\delta,r}^+$ consists of segments 
{of curves either stable or horizontal}.
 Analogously, the singular set of the restriction  $G_{\delta,r}^{-1}: H_{\delta,r}^+ \rightarrow H_{\delta,r}^-$ consists of segments {of curves either unstable or horizontal}.
\label{prop:sing}
\end{proposition}

\begin{proof} 
From Equations  \ref{eqn:T+}, \ref{eqn:T-} or also from general results about the set of singularities of billiards \cite{chernmark}, the singular set consists 
of an union of compact arcs of ${\mathcal C}^\infty$ curves with no other intersection than its endpoints. 
As the curves in $\partial M_{\inn}$ are horizontal, to characterize the singular set of $G$ we only need to analyze curves in $\partial G^{-1}(M_{\inn} )$.
Let us consider a smooth component of  
{$(\mathcal{S}_{\inn}^-\cap H_{\delta,r}^-)\cap \partial G^{-1}(M_{\inn} ) $}.
Any such a curve,  is the pre-image $G^{-1}(\ell)$ of some (horizontal) curve $\ell\subset \partial M_{\inn}$. 
Taking a  {sequence} $\{ \ell_{n} \}$ 
of horizontal curves in $H_{\delta,r}^+ \backslash \mathcal{S}_{\inn}^+$ converging, in $\mathcal{C}^1$ topology, to the curve $\ell$ as $n \rightarrow \infty$, from  Proposition~\ref{prop:us-curves}, $\{G^{-1}(\ell_{n})\}$ is a sequence of stable curves approaching $G^{-1}(\ell)$ as $n \rightarrow \infty$. This  implies that $G^{-1}(\ell)$ is a stable curve. 
\end{proof}

%% file: 4-pontosnormais.tex
\section{Normal Periodic  Points}
\label{sec:normal}

Our strategy to obtain both hyperbolicity and non hyperbolicity is based on the study of the behavior of the first return to the obstacle map $G$ in the neighborhood of some particular periodic orbits known as {\em normal orbits} \cite{foltin}. 
A normal periodic trajectory leaves the obstacle $\alpha$ in the normal direction and, after hitting the exterior circle $\gamma$ at $m+1 $ points ($m\ge1$), 
collides with the obstacle again in the normal direction and, therefore, the same path is traversed with reversed orientation giving rise to an orbit of period 2 for $G$ (or period $2(m+2)$ for $T$)  as shown in Figure~$\ref{fig:normaltraj}$.
{The two  2-periodic trajectories, one from $\omega=\pi$ and the other from $\omega = 0$, which exist for any values of the parameter, are also normal orbits with $m=0$. Both correspond to fixed points of $G$. 
There are also normal periodic orbits with more intermediate hits on the obstacle between the two normal hits, as well as non periodic trajectories with only one normal hit on the obstacle.  
However, we will not consider these two last kind of normal trajectories and, unless specified, we will use the term {\em normal orbits (or trajectories)} only to refer to the two 2-periodic trajectories and to trajectories with exactly two (normal) impacts with the obstacle.

The annular billiard has many normal orbits, and in fact their number increases as $r$ decreases. 
Examples of normal orbits may be constructed in the annular billiard from trajectories leaving  the obstacle in the normal direction and colliding with the external boundary with a rational angle  $\theta= \frac{p}{q}\pi$ .
This situation corresponds to a piece of a trajectory in the circular billiard passing twice through the center of the obstacle.  
Clearly, this construction produces a normal periodic trajectory in the annular billiard for 
{every}
$r$ small enough (Figure \ref{fig:normaltraj}).
It is also clear that the path, and so the period of a normal orbit, for a given $\delta$, remains unchanged as $r$ decreases. As $r\to 0$, each rational $\theta$ will define a normal orbit, implying that the number of normal (periodic) orbits tends to infinity in this limit. 
{The unlimited increase of the number of normal orbits is fundamental in our arguments.}

In an abuse of language, for a given $\delta$, we will call a point $(\omega, 0)$ simply a {\em a normal point} if it corresponds 
to a normal orbit (with two normal hits on the obstacle) for $r$ small enough,
even though strictly speaking, any point $(\omega, 0)$ corresponds to a trajectory leaving the obstacle in the normal direction.

\begin{figure}[h]
\begin{center}
{\includegraphics[width=0.25\hsize]{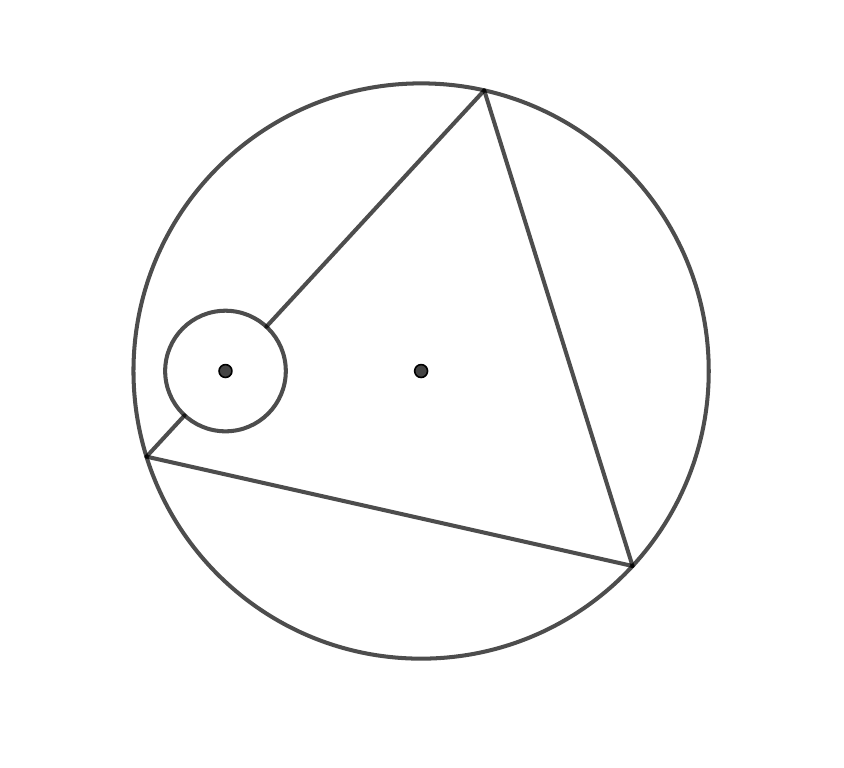}}
\hfill
{\includegraphics[width=0.25\hsize]{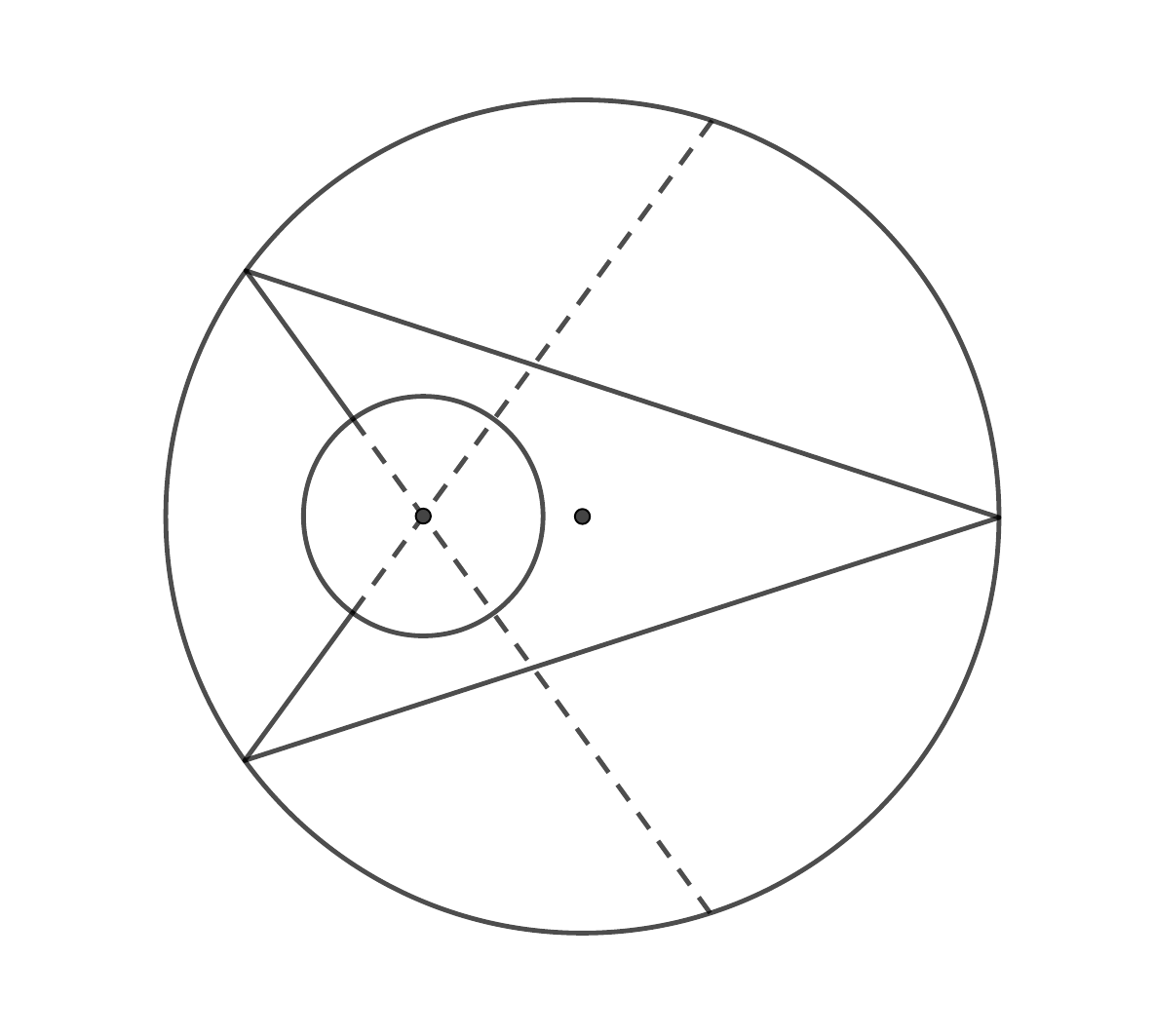}}
\hfill
{\includegraphics[width=0.25\hsize]{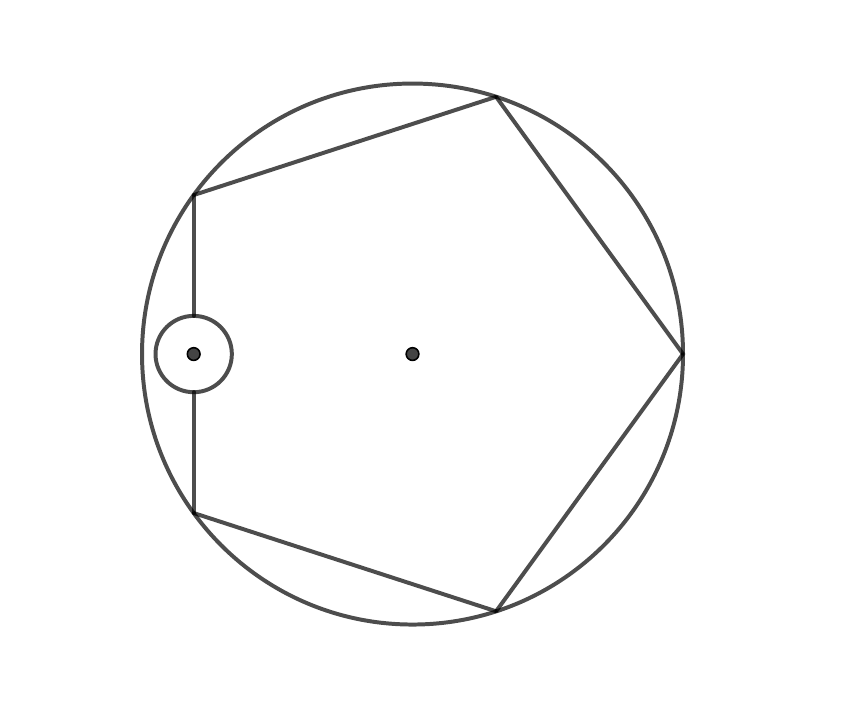}}
\end{center}
\label{fig:normaltraj}
\caption{Normal trajectories}
\end{figure}

In general, as the curve $L^0 \subset  M_{\inn}$ defined by \ref{eqn:L0} denotes the set of orthogonal collisions with the obstacle, 
normal points $(\omega,\beta=0)$ correspond to the intersection $L^0 \cap G^{-1} (L^0) \subset M_{\inn}$.
It follows that $(s,\theta) = T(\omega,0) \in L_{\delta}^+ \cap F^{-m}(L_{\delta}^-) \in M_{\out}$ and so normal points correspond to the
solutions of the system
\begin{equation}
\begin{array}{cl}
 L^+_{\delta}:& \sin\theta   -  \delta \sin \omega =0 \\
 F^{-m}(L_{\delta}^-):&\sin\theta  -  \delta  \sin( \omega - (m+1) (\pi -2\theta))=0
\end{array}
 \label{eqn:pontosnormais}
\end{equation}
{It is worthwhile to notice that normal orbits are symmetric, in the sense that if $(s,\theta)$ belongs to the orbit, so does the point $(s,-\theta)$} 

A normal trajectory of period $2 (m+2)$  is specified by 
$$\left \{ (\omega_0, \beta = 0), (s_0,\theta_0), \ldots,  (s_m,\theta_m), (\omega_1, 0) \right \}
\hbox{ with $s_k = s_0 + k (\pi - 2 \theta_0)$ and $\theta_k = \theta_0$}
$$
with $s_0 = -\theta_0 -\omega_0$, $\omega_1 =   \theta_0 - s_m$
and where $\omega_0$ and $\theta_0$ must satisfy the system \ref{eqn:pontosnormais} above. 
In particular, for any fixed $m \ge 1$, this  system 
has a solution with $\theta$  rational, i.e. a rational multiple of $\pi$.
On the other hand, a solution of the above system  will be a normal orbit in the annular billiard if
$| \sin \theta_k + \delta \sin (\theta_k - s_k) | > r $ for $ 0 \le k < m$, which clearly is verified for any  $r$ small enough. This shows that for any $\delta$ fixed, the number of normal periodic orbits goes to infinity as $r$ decreases to zero.

\begin{definition}
A point $(\omega,0) \in L^0 \cap G^{-1}(L^0)$ is a transverse (resp. tangent) normal point if the intersection is transverse (resp. tangent). 
\end{definition}

{Whether a normal point $(\omega,0)$  is transverse or not depends on the intersection} $T(\omega, 0) \in L_\delta^+ \cap F^{-m}(L_{\delta}^-)$. 
From \ref{eqn:pontosnormais} we have that a tangency occurs if and only if 
\begin{eqnarray}
&&\cos \theta (\cos \omega- \cos (\omega - (m+1)(\pi - 2 \theta)) 
= 2 (m+1) \delta \cos \omega \cos (\omega - (m+1)(\pi - 2 \theta))  \nonumber
\\
&& \hbox{with } \sin \omega = \sin  (\omega - (m+1)(\pi - 2 \theta)) \label{eqn:transverso}
\end{eqnarray}
This implies that tangencies are given by
\begin{equation}
\cos \omega = 0 \  \ \hbox{\  or  } \  \ \frac{\delta \cos \omega}{\cos \theta} = \frac{-1}{m+1}
 \label{eqn:cond-tangent}
\end{equation}
It follows from Lemma~\ref{lem:zeta0} that, for $(\delta,r) \in \Omega_{*}$, all normal points in $H^-_{\delta,r}$ 
are transverse.
Outside $H^-_{\delta,r} \cup H^{+}_{\delta,r}$, we will consider only tangent normal points given by $\cos \omega = 0$, as on Figure~\ref{fig:normaltraj} (R). As already observed in Section~\ref{sec:cones} this last condition implies that $a_{21} = 0$ in the tangent map $DG$ and so these trajectories represent an obstruction to hyperbolicity as obtained there.
In fact, we shall see that  transverse normal points give rise to hyperbolicity, while tangent normal points are related to 
non-hyperbolic dynamics.

We emphasize that trajectories of normal points depend only on $\delta$ (continuously) and do not depend on $r$, since this variable does not intervene in the system  \ref{eqn:pontosnormais} or  
equations \ref{eqn:transverso} and \ref{eqn:cond-tangent}. This is an important remark, as it allows us to use the limit $r \to 0$.

%% file: 5-hyperbolic.tex
\section{Hyperbolic Sets around Transverse Normal Points}
\label{sec:hyperbolic}

In this section we will prove Theorem $\ref{thmA}$ by exhibiting a set of parameters $\Omega_0$ where each first return to the obstacle map  $G_{\delta,r}$ has a horseshoe $\Lambda_{\delta,r}$. The point $(\delta=1,r=0)$ is an accumulation point of the set $\Omega_0$ and  
{the family of horseshoes} $\Lambda_{\delta,r}$ converge to the entire phase space as $(\delta,r) \rightarrow (1,0)$. 
The construction of the horseshoes follows standard arguments as in \cite{hs} and uses,  besides the preservation of cones, the geometric properties of the maps in the neighborhood of transverse normal points, which we describe bellow.

Following the construction in Section \ref{sec:normal}, given $\delta  \in [1/\sqrt{2}, 1)$ we can choose $\omega_i $, such that $(\omega_i,0)$ is a transverse normal point for every $r$ sufficiently small .
We denote by $S_i$ the closure of the connected component of $M_{\inn} \backslash \mathcal{S}_{\inn}^-$ containing $(\omega_i,0)$ 
so, all the points in  ${ int }  (S_i)$ have the same returning time $\nu = \nu(\omega_i,0)$.
The normal trajectory of $(\omega_i,0)$ is 
$2\nu$ periodic and has only two collisions with the obstacle. Since by definition, normal trajectories have no tangential collisions with the obstacle,  the billiard map $T$, and so the first return map $G$, is a ${\mathcal C}^{\infty}$ diffeomorphism in  ${ int }  (S_i)$.
Again, we often omit the dependence on $\delta$ and $r$ of the maps and sets, however we stress that most of the properties of normal transverse periodic orbits depend only on $\delta$  and are actually continuous on this parameter.
A key point in our geometric construction of horseshoes is that, for small $r$,  $S_i$ and $U_i=G(S_i)$ are essentially parallelograms
with two sides in the distinct components of $\partial M_{\inn}$.
This geometric concept will be important in our arguments.
\begin{definition}
A compact connected set {$ \mathbf{S} \subset M_{\inn}$}  is {\em essentially a parallelogram} if its boundary 
is the union of four distinct curves that are $\mathcal{C}^1$ close to the sides of a parallelogram. 
\end{definition}
\begin{definition}
A compact connected set
 $ \mathbf{S} \subset   M_{\inn}$  bounded by two  disjoint  stable (unstable) curves connecting the two opposite components of $\partial M_{\inn}$, will called a {\em stable (unstable) strip}. 
\end{definition}
The expression {\em connecting} $\partial M_{\inn}$ will be always mean {\em connecting the two different components of} $\partial M_{\inn}$.

\begin{lemma}  For each fixed eccentricity $\delta  \in [1/\sqrt{2}, 1)$, and any normal point  $(\omega_i,0)$, with  
$|\sin \omega_i| <   \delta$,
there is  $r_i$
such that for all $r \leq r_i$ the following properties hold
\begin{enumerate} 
\item $S_{i} \ni (\omega_i,0)$ is a stable strip  bounded by two stable curves in $G^{-1}(\partial M_{\inn})$ connecting
the two distinct components $\partial M_{\inn}$ and converging in the ${\mathcal C}^1$ topology to the  decreasing straight line $J_i^-= \{(\omega,\beta): \omega + \beta= \omega_i\}$ as $r\rightarrow 0$.
\item  $U_{i} \ni G(\omega_i,0) = (\hat \omega_i,0)$ is an unstable strip bounded by two unstable curves in $G(\partial M_{\inn})$ connecting  $\partial M_{\inn}$ and converging in the ${\mathcal C}^1$ topology to the increasing straight line $J_i^+= \{(\omega,\beta): \omega - \beta= \hat \omega_i\}$ as $r\rightarrow 0$.
\end{enumerate}
\label{lem:strips}
\end{lemma}

\begin{proof} By reversibility, it is enough to prove $1$. 
The  first return to the obstacle time of $(\omega_i,0)$ is  $\nu(\omega_i,0) = m_i+2$ for some integer $m_i \ge 0$.
So, the  first return map restriction $G: S_i \rightarrow U_i$ decomposes as  $G=T \circ F^{m_i} \circ T$ and by definition, $S_i$ is the connected component containing $(\omega_i,0)$ of the set
$
M_{\inn} \cap G^{-1}(M_{\inn}) \backslash (T^{-2}(M_{\inn}) \cup ... \cup T^{-m_i-1}(M_{\inn}))
$

{
To describe $S_i$, we consider its image  $T(S_i)$ which is a subset of $M_{\inn}^+ \cap F^{-m_i}(M_{\inn}^-)$. If $V_i$ denotes the connected  component of $M_{\inn}^+ \cap F^{-m_i}(M_{\inn}^-)$ containing $T(\omega_i,0)$,
{it is clear that $T(S_i) \subset V_i$} and we have
\begin{equation}
T(S_i) = V_i \backslash (F^{-1}(M_{\inn}^-) \cup \ldots F^{-m_i+1}(M_{\inn}^-)) 
\subset M_{\inn}^+ \cap F^{-m_i}(M_{\inn}^-)  \subset M_{\out}
 \label{eqn:Vx}
\end{equation}

As observed in the end of Section \ref{sec:preliminares}, 
$\partial M_{\inn}^{\pm} \underset{r \to 0}{\longrightarrow}  L_{\delta}^{\pm}$, so for $j= 0,\ldots,m_i$
we have 
\begin{eqnarray}
&& F^{j} (\partial M_{\inn}^+) \underset{r \to 0}{\longrightarrow} F^{j}( L_{\delta}^+)
\ , \
F^{-j} (\partial M_{\inn}^-) \underset{r \to 0}{\longrightarrow} F^{-j}( L_{\delta}^-)
 \hbox{ \ \ in $C^\infty$ topology}
\label{eqn:ConvergeCurves}
\\
&& F^{\pm j} (M_{\inn}^\pm) \underset{r \to 0}{\longrightarrow} F^{\pm j}( L_{\delta}^\pm)
 \hbox{ \ \  in the Hausdorff set distance} 
\nonumber
\end{eqnarray}

The set $H$, defined in Section \ref{sec:cones}, is the horizontal strip $|\sin \theta| < \delta^2$, so the choice $|\sin \omega_i| $ means that 
$(\omega_i,0) \in H^{-}$
for any $r$
and so it is a transversal normal point.
It follows that the intersection $L^+_{\delta} \cap  F^{-m_i}(L^-_{\delta})$ at $T(\omega_i,0)$ is also transversal. 
By definition $V_i \subset M_{\inn}^+ \cap F^{-m_i}(M_{\inn}^-)$ so this transversality and  \ref{eqn:ConvergeCurves} imply that
\begin{equation}
F^{j} (V_i) \underset{r \to 0}{\longrightarrow} T^{j+1}(\omega_i,0) 
\hbox{ for $j=0,\ldots,m_i$} 
\label{eqn:ConvergeSets}
\end{equation}
On the other hand, as the first returning time is $m_i+2$, $T^{{m_i}+1}(\omega_i,0) \in M_{\inn}^-$, but $T^j(\omega_i, 0) \notin M_{\inn}^-$ for $j= 1, \ldots,m_i$.  From 
\ref{eqn:Vx} and \ref{eqn:ConvergeSets}, we can choose $r_i$ small enough in order  that if $r \leq  r_i$  we have $F^j(V_i) \cap M_{\inn}^-= \emptyset$ for $j= 1, \ldots , m_i$, {implying that, in fact, $T(S_i)= V_i$}.

Furthermore, the convergence of $T(S_i)$ to $T(\omega_i,0)$  in \ref{eqn:ConvergeSets} and the transversality of the intersection between 
$L_\delta^+$ and  $ F^{-m_i}(L_{\delta}^-) $ at this point, together with 
the convergence of
$\partial M_{\inn}^+$ to  $L_{\delta}^+$  and of $F^{-m_i}(M_{\inn}^-)$ to $F^{-m_i}(L_{\delta}^-) $
in \ref{eqn:ConvergeCurves},
imply that for $r$ small enough $T(S_i)$ is essentially a  parallelogram  bounded by two curves in $\partial M_{\inn}^+$ and two curves in $F^{-m_i}(\partial M_{\inn}^-)$. 
Hence $S_i \subset M_{\inn}$ is a strip  bounded by two curves  
in $T^{-1} \circ F^{-m_i}(\partial M_{\inn}^-) \subset G^{-1}(\partial M_{\inn})$ 
connecting $\partial M_{\inn}$.

Moreover, we have that $T(\omega_i,0) = (s,\theta)$ with $|\sin \theta| < \delta^2$, which means that $T(\omega_i,0) \in H$ and clearly we can set $r_i$ such that $T(S_i)\subset H = T(H^-)$ implying  that
$S_i \subset H^-$.
 Thus the two curves connecting  $\partial M_{\inn} \subset \partial S_i $
are stable, since they belong to the singular set  $G^{-1}(\partial M_{\inn}) \cap H^{-}$
 of $G$ as discussed in Proposition \ref{prop:sing}.
This proves that $S_i$ is a stable strip.
}

To prove that $\partial S_i \to J^-_i$, we refer to Proposition  \ref{prop:us-curves}. 
The two opposite stable curves of 
$\partial S_i \subset S_i \cap G^{-1}(\partial M_{\inn})$ converge, in the ${\mathcal C}^1$ topology, to straight lines of slope $-1$ as $r\to 0$.
{Now, consider an horizontal segment $\ell_\beta$  connecting these two curves. 
Its image $G(\ell_{\beta}) \subset U_i$ connects the components of $\partial  M_{\inn}$ and $|G(\ell_{\beta})| > \rho |\ell_{\beta}|$.
On the other hand, as  $G(\ell_{\beta})$ is a Lipschitz curve with a constant close to 1 connecting the boundaries,  its length is less than some constant close to $\sqrt{2}$ and so for any $\beta$, $|\ell_{\beta}| \lesssim \frac{\sqrt{2}}{\rho}$.
It follows that $|\ell_{\beta}|  \to 0 $ uniformly as ${r \to 0}$, implying that $S_i \rightarrow J_{i}^{-} $ as $ r \rightarrow 0$.}
\end{proof}

The construction of the horseshoes for $G$ is based on sets of 
transverse normal points  $X_{\delta}$ 
which we describe bellow.

From  the results of Section \ref{sec:normal}, for any arbitrary $\delta  \in [1/\sqrt{2}, 1)$,  there is a dense set
of points in 
$H^-_{\delta,r} \cap L^0 = \{ (\omega,0): |\sin \omega| < \delta \}$ 
that will give rise to transverse normal orbits as $r \rightarrow 0$.
In this dense set of transverse normal points, we can choose {\ver a set with $n_{\delta}$} points
such that the ${\mathcal S}^1$-distance between any two adjacent points is less than $d_{\delta} < \pi - 2 \arcsin \delta$
\begin{equation}
 X_{\delta}= \left \{ (\omega_1,0), (\omega_2,0), \ldots (\omega_{n_{\delta}},0)  \right \} 
\subset H^-_{\delta,r}  
\label{eqn:X}
\end{equation}
We observe that $ \pi - 2 \arcsin {\delta}$ is the length of each of the two disjoint components of the complement $L^0 - {H^-_{\delta,r}}$
which are located around $(\omega = \pm \pi/2,0)$.
By including images, we can assume that $X_{\delta}$ is invariant under $G$.
With this choice, the invariant set $X_{\delta}$ becomes dense in $L^0$ as $\delta \to 1$ 
{\ver and obviously  $n_{\delta} = \# (X_{\delta}) \to \infty$}.
It is important to notice that the set $X_{\delta}$ is robust on $\delta$ and does not depend on $r$ as long it is sufficient small. In particular $H^-_{\delta,r} \cap L^0$ does not depend on $r$.

\begin{lemma}
\label{lem:horseshoe}
For each fixed $\delta \in [1/ \sqrt{2},1)$ there is $r_{\delta}$ such that for any  $r \in (0, r_{\delta}]$ the  map $G_{\delta,r} $ 
has a 
{locally maximal transitive hyperbolic} set $\Lambda_{\delta,r}$ such that
\begin{enumerate}
\item The restriction $G_{\delta,r}: \Lambda_{\delta,r} \rightarrow \Lambda_{\delta,r}$
 is conjugated to a  sub shift in the space  of sequences of $n_{\delta}$ symbols. 
{\ver Moreover, $n_{\delta} \to \infty$ as $\delta \to 1$.}
\item  The set $\Lambda_{\delta,r}$ is $d_{\delta}$-dense in $M_{\inn}$ 
with $d_{\delta} \rightarrow 0$ as $\delta \rightarrow 1$.
\item For any pair $r' \ne r$ in  $(0,r_{\delta}]$  
the set $\Lambda_{\delta,r'}$ is the hyperbolic continuation of the set  $\Lambda_{\delta,r}$.
\end{enumerate}
\end{lemma}

\begin{proof}
Given two distinct normal points $(\omega_i,0)$ { and } $(\omega_j,0)$  in $X_{\delta}$,  let $(\hat \omega_i, 0)=G(\omega_i,0) \in X_{\delta}$. 
We will investigate the  intersection $U_{i} \cap S_{j}$ where $U_i \ni (\hat \omega_i,0) $ and $S_j \ni (\omega_j,0)$.
It follows from Lemma~\ref{lem:strips},  that for small $r$  this intersection is related to the intersection of the
lines $J_i^+ \ni (\hat \omega_i,0) $ and  $J_j^- \ni (\omega_j,0)$. 
It is obvious that $J_i^+\cap J_j^-$ consists of a single point in the interior of $M_{\inn}$, unless $\hat \omega_i = \omega_j  + \pi$ in which case it consists of two points in the distinct components of $\partial M_{\inn}$. 
As a consequence, if $ \hat\omega_i - \omega_j \neq  \pi $, for $r$ small enough, $U_i \cap S_j$ is essentially a parallelogram bounded by two unstable curves in $\partial U_{i}$ and two stable curves in $\partial S_{j}$.

In what follows we consider  $\displaystyle  0<  r_{\delta} < \min_{i = 1 \ldots n_{\delta}} r_i$, where $r_i$ is given by Lemma~\ref{lem:strips}. 
Clearly we can also assume that $r_{\delta}$ is small enough so that $U_i$ effectively crosses
 $S_{j}$  whenever $ \hat \omega_i - \omega_j \neq \pi$. 
{It is also clear that by at least three points in the set $X_{\delta}$ we will always have such crossings.}

Given the set $X_{\delta}$ of $n_{\delta}$ transverse normal points, let us consider
$ \Sigma=\{1,...,n_{\delta}\}^\mathbb{Z}$,  the space of sequences $a=\{a_i\} _{i \in \mathbb{Z}}$ of 
$n_{\delta}$ symbols and the shift map $\sigma$ on it.
We define
{$\widetilde{\Sigma} \subset \Sigma$}, as the $\sigma$-invariant subset of sequences  such that, for any $i \in \mathbb{Z}$,  $U_{a_i}=G(S_{a_{i}})$ crosses $S_{a_{i+1}}$. 
Now, given a sequence $a  = \{a_i\}_{i \in \mathbb{Z}} \in \widetilde{\Sigma}$ we define the sets
\begin{equation}
S_{a}^n= \bigcap_{j=0}^nG^{-j}(S_{a_{j}}) 
\hbox{ \ and \  }  
U_{a}^n= \bigcap_{j=0}^n G^{j}(U_{a_{\text{-}{(j+1)}}})
\hbox{ \ , for $n \geq 0 $} 
\label{eqn:SnUn}
\end{equation}
So,  $S_{a}^{n+1} \subset S_{a}^n \ldots \subset S_a^0 = S_{a_0}$ and 
$U_{a}^{n+1} \subset U_{a}^n \ldots \subset U_a^0 = U_{a_{\text{-1}}}$.

We will show that each $S_{a}^n$ is 
a {stable} strip
bounded by stable curves in $G^{-(n+1)}(\partial M_{\inn})$ {connecting $\partial M_{\inn}$} and similarly that each $U_{a}^n$ is an {unstable} strip bounded by unstable curves in $G^{(n+1)}(\partial M_{\inn})$
{also connecting $\partial M_{\inn}$}.
This is obviously true for $S^0_a = S_{a_0}$ and
$U^0_a = U_{a_{\text{-1}}}$ by Lemma \ref{lem:strips}. Moreover, we note that by definition  $S^0_a$ crosses $U^0_a$ and so,
when  all  the strips are stable or unstable, 
the intersections $S^{n}_a$ and $U^{n}_a$ also cross.

We have that 
\begin{equation}
S_{a}^1= S_{a_0} \cap G^{-1}(S_{a_1}) = G^{-1}(G(S_{a_0}) \cap S_{a_1}) = G^{-1}(U_{a_0} \cap S_{a_1})
\label{eq:S1a}
\end{equation}
As $U_{a_0}$ is a strip bounded by two {unstable} curves in $G(\partial M_{\inn})$ and  
$S_{a_{1}}$ is a strip bounded by  two {stable} curves in $G^{-1}(\partial M_{\inn})$,
$U_{a_0} \cap S_{a_{1}}$ is essentially a parallelogram  bounded by two curves in  
$\partial U_{a_0}\cap G(\partial M_{\inn})$ and two curves in $\partial S_{a_1} \cap G^{-1}(\partial M_{\inn})$.
Its image under $G^{-1}$ is also essentially a parallelogram, bounded by two curves in distinct components of 
$\partial M_{\inn}$ and two stable curves in $G^{-2}(\partial M_{\inn})$. Hence $S_{a}^1$ is a stable strip in $S^0_{a} = S_{a_0}$ bounded by two opposite curves in $G^{-2}(\partial M_{\inn})$.
A similar argument shows that $U_{a}^1$ is an unstable strip in $U^0_{a} = U_{a_{\text{-1}}}$ bounded by two opposite curves in $G^{2}(\partial M_{\inn})$.

The same construction can be applied to $S^n_a$ for  $n>1$, and also for $U^n_a$. For instance,
$$S^2_a= S_{a_0} \cap G^{-1}(S_{a_1}) \cap G^{-2}( S_{a_2}) = S^1_a \cap  G^{-2} (S_{a_2})
= G^{-2}(G^2(  S^1_a) \cap  S_{a_2})
$$
By equation \ref{eq:S1a} we have $G^2(S_a^1)=G (U_{a_0} \cap S_{a_{1}}) $.
As $S_a^1$ has two boundaries in $G^{-2}(\partial M_{\inn})$,  $G^2(S_a^1)$ is a strip and, as the two other boundaries of $S_a^1$ are in $\partial M_{\inn}$, their image under $G^{2}$ are unstable curves.  So $G^2(S_a^1)$ is an unstable strip in $G(S_{a_1}) =U_{a_1}$  and it must cross  $S_{a_2}$ implying that the intersection $G^2(S_a^1) \cap S_{a_2}$ is also essentially a parallelogram with two boundaries in $G^2(\partial M_{\inn})$ and two in $G^{-1}(\partial M_{\inn})$.
It follows that $S_{a}^2 \subset S_{a}^1$ is a strip bounded by two stable curves in $G^{-3}(\partial M_{\inn})$. 

By induction we assume that $S_{a}^{n-1}$ is a (stable) strip bounded by stable curves in $G^{-n}(\partial M_{\inn})$.
The definition \ref{eqn:SnUn} can be written as
\begin{eqnarray*}
S^{n}_a &=& S_{a_0} \cap G^{-1}(S_{a_1}) \ldots \cap G^{-(n-2)}(S_{a_{n-2}})   \cap G^{-(n-1)}(S_{a_{n\text{-1}}}) \cap G^{-n}(S_{a_n}) \\
&=&   S^{n-1}_a \cap G^{-n}(S_{a_n} ) =  G^{-n}  (G^{n}(S^{n-1}_a) \cap S_{a_n} ) 
\end{eqnarray*}
The induction hypothesis implies that  ${G^n}(S^{n-1}_a)$ is a strip bounded by two unstable curves in $G^n(\partial M_{\inn})$ connecting $\partial M_\alpha$. 
As by definition 
$S^{n-1}_a \subset G^{n-1}(S_{a_{n\text{-1}}})$, we have that $G^n(S^{n-1}_a) \subset G(S_{a_{n\text{-1}}}) = U_{a_{n\text{-1}}}$. It follows that 
$G^{n}(S^{n-1}_a) \cap S_{a_n} $ is essentially a parallelogram with two boundaries in $G^{n}(\partial M_{\inn})$ and two in 
$G^{-1}(\partial M_{\inn})$ and so taking its image  
under $G^{-n}$, we obtain that $S^{n-1}_a \cap G^{-n}(S_{a_n})  = S^{n}_a$  is a stable strip with boundaries in $G^{-(n+1)} \partial M_{\inn}$.

\begin{figure}
\begin{center}
\includegraphics[width=0.9\hsize]{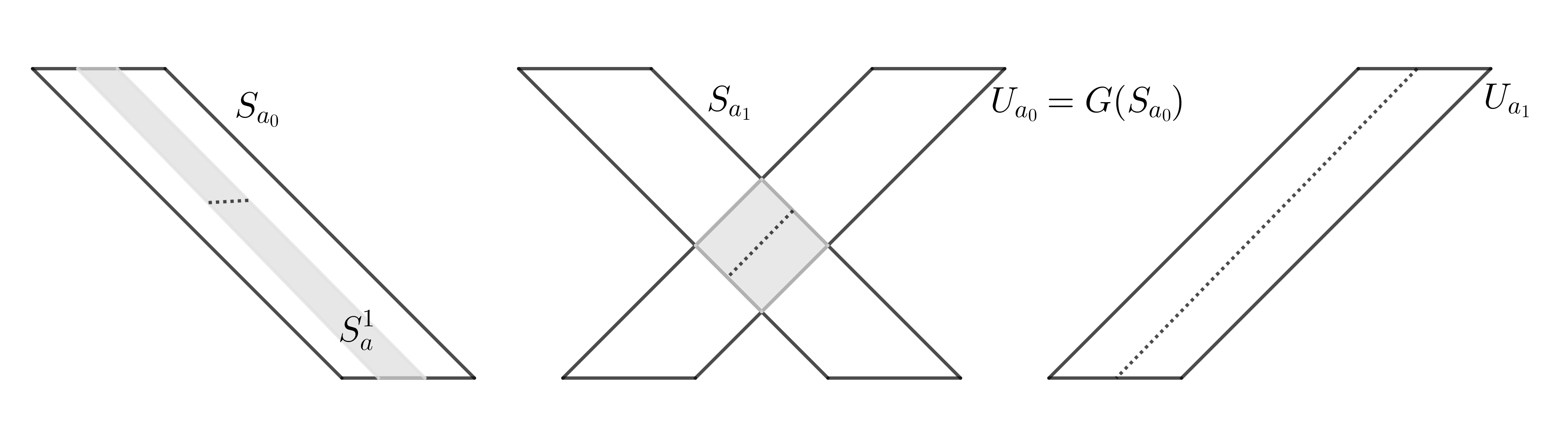}
\end{center}
\caption{Construction of the sets $S_n$}
\end{figure}


Let us consider $S^n_a$ and an horizontal segment $\ell_{\beta}$ connecting the two opposite stable curves of 
$\partial S^n_a \subset G^{-(n+1)}(\partial M_{\inn})$.
Using an argument similar to the one 
at the end of the proof of Lemma \ref{lem:strips},  we have that the horizontal width of the strip 
$S^n_a$ is bounded by $|\ell_{\beta}| \lesssim \frac{\sqrt{2}}{\rho^n} \to 0 $ as $n \to \infty$.
This  convergence together with the properties of stable curves already stated imply that,
for any  $a \in \tilde \Sigma$,
$S^{\infty}_a = {\bigcap_{n=0}^{\infty}} S^n_a$ is a decreasing $1/c_1$-Lipschitz curve connecting $\partial M_{\inn}$. From reversibility $U_{a}^\infty$ is an increasing $1/c_1$-Lipschitz curve connecting $\partial M_{\inn}$.
Hence, to each $a \in \widetilde{\Sigma}$ corresponds a unique point $S_a^\infty \cap U_a^\infty$ in $M_{\inn}$ and we can define a map  $h: \widetilde{\Sigma} \rightarrow M_{\inn}$ given by $h(a)= S_{a}^\infty \cap U_{a}^\infty$.  Standard arguments \cite{hs} show that $h$ is a homeomorphism onto its image. 

To obtain the hyperbolic set, we define $\Lambda_{\delta,r}=h(\tilde \Sigma)$ which is a compact $G$-invariant set in $H_{\delta,r}^-$. The preservation of cones in  $H_{\delta,r}^-$  (Corollary \ref{Coro:Cones})
implies that $\Lambda_{\delta,r}$ is a hyperbolic set for $G$.  Moreover, the definition of $h$  implies that $G$ restricted to $\Lambda_{\delta,r}$ is conjugated to the shift map $\sigma: \tilde \Sigma \rightarrow \tilde \Sigma$.

For small $r$, the sets $S_i$ and $U_i$ are respectively close to the lines $J^-_i$ and $J^+_i$.  
The points $J^-_i \pitchfork J^+_i = (\omega_i,0) \in X_{\delta}$  are $d_{\delta}$-dense in $L^0$ 
and so we have a square lattice of lines $J^-_i$ and $J^+_k$ with $i, \,k=1\ldots n_{\delta}$, which nodes 
$J^-_i \pitchfork J^+_k$ are $\frac{d_{\delta}}{\sqrt{2}}$-dense in $M_{\inn}$.
By definition, the points in $\Lambda_{\delta,r}$ are close to the nodes
in the interior of $M_{\inn}$  
and 
therefore we can set $r_{\delta}$ such that the hyperbolic set itself is $d_{\delta}$-dense  
in $M_{\inn}$.
It is clear from the construction that, for each fixed $\delta$ and each $r \in (0,r_{\delta}] $,
$\Lambda_{\delta,r}$
is a locally maximal transitive hyperbolic set.
This implies that it has a continuation in $r$ which in turn is a locally hyperbolic.
{More precisely, there is an open set $V \ni r$, 
such that for any $G_{\delta,r'}$ with $r' \in V$,
the hyperbolic set  $\Lambda_{\delta,r'}$ is the continuation of  $\Lambda_{\delta,r}$.} 
As the argument holds for any $r$,  we can take $ V =  (0,r_{\delta}]$. This proves item 3.
\end{proof}

\begin{lemma} 
\label{lem:extendedhorseshoe}
For any  $\delta \in [1/\sqrt{2},1)$ there is a set of parameters $R_{\delta}=(\delta-\epsilon_{\delta},\delta+ \epsilon_{\delta}) \times (0,r_{\delta}] \subset \Omega_*$ 
{such that} for any  $(\delta,r) \in R_{\delta}$ the map $G_{\delta,r}$ 
has
a locally maximal transitive hyperbolic set
$\Lambda_{\delta,r}$. 
Moreover, if  $(\delta',r')$ and $(\delta'', r'')$ are in 
$R_{\delta}$, the set $\Lambda_{\delta'',r''}$ is the {continuation} of $\Lambda_{\delta',r'}$.
\end{lemma}

\begin{proof} 
Given $\delta$, the continuity and the transversality imply that if $\tilde \delta \approx \delta$, for each $(\omega_i,0) \in X_{\delta}$ we can find $\tilde\omega_i \approx \omega_i $ (called the continuation of $\omega_i$)  such that $(\tilde\omega_i, 0 )$ is also a transverse normal point of $G_{\tilde \delta, r}$.
In fact,  there is  $\epsilon_{\delta}$ and we can adjust $r_{\delta}$ 
such that if $|\widetilde{\delta}-\delta| {<} \epsilon_{\delta}$ then 
the set  $X_{\tilde \delta}$ obtained by continuation of $X_\delta$ contains also $n_{\delta}$ points which are normal transverse for any $0 < r\leq  r_{\delta}$.  
If necessary, we can take smaller $r_{\delta}$ and $\epsilon_{\delta}$, in order to ensure that  
$R_{\delta} = (\delta-\epsilon_{\delta},\delta+ \epsilon_{\delta}) \times (0,r_{\delta}] \subset \Omega_*$ (as given by equation \ref{eqn:omega*}).  

Applying Lemma \ref{lem:horseshoe} we construct a locally maximal hyperbolic set $\Lambda_{\bar \delta,r}$ for any $r  \leq r_{\delta}$. 
It is obvious that if $\delta'$ and $\delta''$ are close, the square lattices obtained from $X_{\delta'}$ and $X_{\delta''}$, as in Lemma \ref{lem:horseshoe} are close for small value of $r$ and so are the
hyperbolic sets
 $\Lambda_{\delta',r}$ and $\Lambda_{\delta'',r}$.
 It follows from the continuity on $\delta$ and $r$
 and the uniqueness of the hyperbolic continuation that for any $(\delta',r')$ and $(\delta'', r'')$ in $R_{\delta}$, the set $\Lambda_{\delta'',r''}$ is the continuation of $\Lambda_{\delta',r'}$.
\end{proof}


We now proceed to the 
\begin{proof}[Proof of Theorem $\ref{thmA}$]
We can choose a {countable} covering of $( \frac{1}{\sqrt{2}},1)$ by 
intervals  $(\delta_k-\epsilon_{\delta_k},\delta_k+ \epsilon_{\delta_k})$ 
and take
$R_{\delta_k}$ as defined  in  Lemma~\ref{lem:extendedhorseshoe} and where
$\{\delta_k \}_0^{\infty}$  is a strictly increasing sequence 
with $\delta_0 > 1/\sqrt{2}$ and $\delta_k$ converging to 1 as $k \to \infty$.

To any $\delta \in (\frac{1}{\sqrt{2}},1) $, we assign $R_{\delta_k}$ by taking  $k$ as  the smallest integer such that 
$|\delta-\delta_{k}| { <} \epsilon{_k}$.
By Lemma~\ref{lem:extendedhorseshoe}, for each 
$(\delta,r)  \in \Omega_0=\bigcup_{k \geq 0} R_{\delta_k}  \subset \Omega_*$  the map $G$ has a hyperbolic set $\Lambda_{\delta,r}$. The set $ \Omega_0$ is a connected set with non empty interior such that the point $(\delta=1,r=0)$ is in its boundary.

With the choice above, given any $(\delta,r) \in R_{\delta_k} \subset \Omega_0$, we have a piecewise continuous family of horseshoes $\Lambda_{\delta,r}$ where the map $G$
is conjugated to a shift of $n_{\delta_k}$ symbols, with $n_{\delta_k} \to \infty$ .
Each horseshoe is $d_{\delta_{k}}$-dense in $M_\alpha$ with $d_{\delta_{k}} \to 0$, 
as $\delta_k \to 1$. 
\end{proof}

{\ver At this point it is important to notice that a different choice of the set $X_{\delta}$ yields in principle to a different hyperbolic set and so we in fact could up with many of them.
}

We conclude this section with the description of the stable and unstable manifolds of the hyperbolic set.
The existence of the singularities of the map $G_{\delta,r}$ implies that the global invariant manifolds of points in 
$\Lambda_{\delta,r} $ are disconnected. In what follows we describe the properties of the connected local manifolds which will be essential in some of our geometric arguments.

Fixing $\delta_0 \in (1/\sqrt{2},1]$ let us consider, as in Lemma \ref{lem:extendedhorseshoe}, the two parameter family of maps $G_{\delta,r}$, $(\delta,r) \in R_{\delta_0}=(\delta_0-\epsilon_0,\delta_0 + \epsilon_0) \times  (0,r_{\delta_0}]$ and a corresponding  two parameter family of hyperbolic sets $\Lambda_{\delta,r}$.

The {\em local  stable invariant manifold} of a point  $z \in \Lambda_{\delta,r}$, denoted by $W_{\rm loc}^s(z)$, is defined as the connected component of the stable manifold of $z$ containing this point. 
It is clear from the proof of Lemma \ref{lem:horseshoe} that it is a ${\mathcal C}^\infty$ stable curve connecting the two different components of the boundary $\partial M_{\inn}$. 
Likewise, the \textit{local unstable manifold} of  $z \in \Lambda_{\delta,r}$, denoted by  $W_{\rm loc}^u(z)$,
is  a ${\mathcal C}^\infty$ unstable curve connecting the different components of the boundary $\partial M_{\inn}$.

Now, still for $(\delta,r) \in R_{\delta_0}$, we can consider a two parameter family of points $z_{\delta,r} \in \Lambda_{\delta,r}$ such that any two points in this family are the continuation of each other.
We refer to such a family as a {\em continuous family}. 
The set of admissible sequences $\tilde \Sigma$ does not depend on the parameters $\delta$ and $r$ as long as they stay in $R_{\delta_0}$ and it is clear from the construction, that the points of a continuous family share the same symbolic representation $a$ in $\tilde \Sigma$.
So, given a sequence $a = (a_0, a_1, \ldots) \in \tilde \Sigma$, we consider the corresponding two parameter family of local stable manifolds $W_{\rm loc}^s(z_{\delta,r})$ related to the associated  two parameter family of points $z_{\delta,r}$. 
For each $(\delta,r)$, $W_{\rm loc}^s (z_{\delta,r})$ belongs to the strip $S_{a_0}$ containing a normal point $(\omega^{a_0}_{\delta},0)\in X_{\delta}$. 
Then, for $\delta$ fixed, as $r \rightarrow 0$
the stable boundary of $S_{a_0}$,
and so the local stable manifold $W_{\rm loc}^s(z_{\delta,r})$,
converges in ${\mathcal C}^1$ topology  to 
the  straight line $\omega + \beta = \omega^{a_0}_{\delta}$.

Thus, for $(\delta,r) \in R_{\delta_0}$ and any fixed  $a \in\tilde \Sigma$  we have a two parameter family of local stable manifolds $W^s_{\rm loc}(z_{\delta,r})$ converging as $(\delta,r) \to (\delta_0,0)$ to the decreasing line  $\omega + \beta = {\omega^{a_0}_{\delta_0}}$. 
Correspondingly, the two parameter family $W^u_{\rm loc}(z_{\delta,r})$ of local unstable manifolds converges to the increasing line $\omega-\beta= { \omega^{a_0}_{\delta_0}}$, as $(\delta,r) \rightarrow (\delta_0, 0)$.

%% file: 6-tangencias.tex
\section{Conservative Newhouse Phenomenon}
\label{sec:tangencias}

In this section we prove Theorems $\ref{thmB}$ and $\ref{thmC}$ by exhibiting a set of parameters
in $\Omega_0$ and  accumulating $(\delta=1,r=0)$ such that the first return to the obstacle map  presents quadratic homoclinic tangencies that unfold generically as the radius of the obstacle  varies. Accumulating this parameter set with  homoclinic tangencies, we find another set where the map has elliptical islands  filling in the phase space as $(\delta,r) \rightarrow (1,0)$.

These phenomena {originate from the bifurcation of}
tangent normal points, defined in Section \ref{sec:normal}. 
As pointed there,
for $ \delta= \sin \frac{p}{q}\pi \equiv \delta_0$,
where $ 0<\frac{p}{q}<1$ is any  rational number,
the point $(\omega_0=-\frac{\pi}{2},\beta_0 = 0) \in M_{\inn}$ is a tangent normal point for $G$.

\begin{figure}[h]
\begin{center}
\includegraphics[width=0.3\hsize]{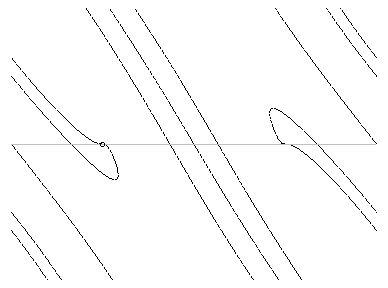}\hskip 1cm\includegraphics[width=0.3\hsize]{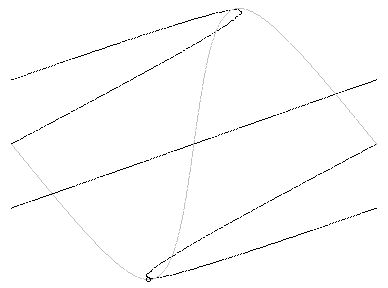}\hfill
\end{center}
\caption{Tangent normal point 
($L^0, G^{-1}(L^0) \subset M_{\inn}$ and $ L^+_{\delta_0}, F^{-m}(L^-_{\delta_0}) \subset M_{\out}$) 
}
\label{fig:normaltangente}
\end{figure}

The local study of the (tangent) intersection
 between the horizontal line $L^0$ and its preimage $G^{-1}(L^0)$ reveals that this  tangency 
is cubic (Figure \ref{fig:normaltangente}) 
and unfolds generically as $\delta$ varies.

For values of $\delta \approx \delta_0$ and  small $r$
the curve $G^{-1}(L^0) $ is  ${\mathcal C}^{ 1}$  close to a segment of the stable manifold in the hyperbolic set 
$\Lambda_{\delta,r}$,
in the neighborhood of the point $(-\frac{\pi}{2}, 0)$.
The local geometric properties of the stable manifold, inherited from the proximity of the tangent normal point, give rise to a quadratic homoclinic tangency, which unfolds generically as r varies. 
The bifurcation of the homoclinic  tangency, implies the appearance of elliptical islands for nearby parameter values\cite{duarte}.

{In the lemmas leading to  the proof of the two theorems we will focus on the neighborhood of the orbit of the tangent normal point and consider two parameter families of maps $G_{\delta,r}$ with  
$(\delta,r) \in R_{\delta_0}$ for different values of $\delta_0$.}
The choice of the set 
$$ R_{\delta_0} = (\delta_0- \epsilon_0,\delta_0 + \epsilon_0) \times (0,r_0]   \subset \Omega_0
$$
as defined in Lemma \ref{lem:extendedhorseshoe},  assures the existence of a continuous family of hyperbolic sets 
$\Lambda_{\delta,r}$.
Eventually we will need to take smaller values of the  constants  $\epsilon_0$ and $r_0$.


{We begin by investigating the bifurcation in $\delta$ of the tangent normal point of some $\delta_0$.}

\begin{lemma} 
\label{lem:BTP}
For $\delta=\delta_0=\sin\frac{p}{q}\pi$, the point $(\omega_0 = -\frac{\pi}{2},\beta_0 = 0) \in L^0 \cap G_{\delta_0,r}^{-1}(L^0)$ is a cubic tangency unfolding into three transverse intersections when $\delta>\delta_0$
and small fixed $r$.
\end{lemma}

\begin{figure}[h]
\includegraphics[width=0.25\hsize]{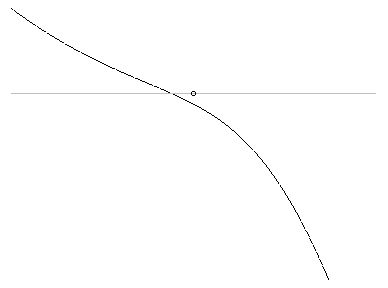}
\hfill
\includegraphics[width=0.25\hsize]{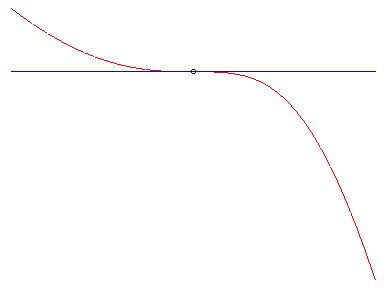}
\hfill
\includegraphics[width=0.25\hsize]{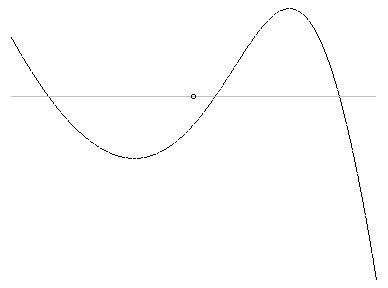}

\includegraphics[width=0.25\hsize]{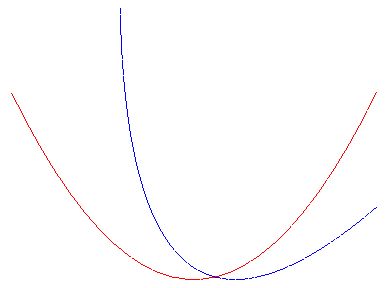}
\hfill
\includegraphics[width=0.25\hsize]{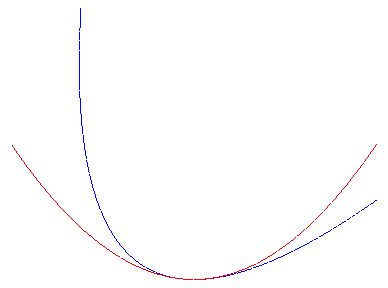}
\hfill
\includegraphics[width=0.25\hsize]{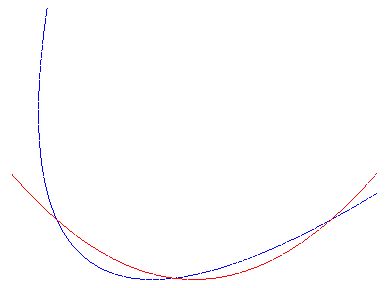}
\caption{The unfolding of the cubic tangency in $M_{\inn}$ (top) and $M_{\out}$ (bottom)}
\label{fig:tangencias}
\end{figure}

\begin{proof}
Assuming that the tangent normal point has return time $\nu(-\frac{\pi}{2}, 0)=m+2$, 
 the first return to the obstacle map $G_{\delta_0,r}$,
{in the connected component  of $M_{\inn}\backslash {\mathcal S}^-_{\inn}$ containing it,} decomposes as
$G_{\delta_0,r}= T\circ F^{m} \circ T$ where $T=T_{\delta_0,r}$. 
If $\delta$ close to $\delta_0$ we have the same return time and so the same decomposition for $G_{\delta,r}$. 
Moreover, in this {neighborhood},  the maps are $\mathcal{C}^\infty$.
 
We will describe the bifurcation of  $L^0 \cap G^{-1}(L^0) $ by looking at the bifurcation of its image 
{$(s_0,\theta_0)$} in
$L_{\delta}^+ \cap F^{-m}(L_{\delta}^-)  \subset M_{\out}$. 
From the definition of  the curves $L_{\delta}^\pm$ and the map  $F$, the points in 
$F^{-m}(L_{\delta}^-) \cap L_{\delta}^+$, correspond to the solutions $(s,\theta)$ of the system  \ref{eqn:pontosnormais} for each $\delta$ 
and any $r$ sufficiently small.
For simplicity we describe the case of odd $m$ (the even case is similar).
Introducing the variable $\varphi = s+\theta = -\omega$ and including the dependence on the parameter $\delta$, the system is written as
\begin{eqnarray}
&L^+_{\delta}: &  A(\varphi,\theta;\delta)  =  \sin\theta +\delta \sin\varphi =0  \label{eqn:A}\\
&F^{-m}(L_{\delta}^-): & B(\varphi,\theta;\delta)  =  \sin\theta +\delta \sin(\varphi-2(m+1)\theta)= 0  \label{eqn:B}
\label{S1}
\end{eqnarray}

{This system has at least a solution $(\varphi_0,{\theta_0};\delta_0)= (\frac{\pi}{2}, -\frac{p}{q}\pi; \sin(\frac{p }{q}\pi))$  which corresponds to the image of the tangent normal point $(\omega_0,\beta_0)=(-\frac{\pi}{2},0)$ since}
\begin{eqnarray*}
 A(\varphi_0,\beta_0;\delta_0)  &=&  \sin(-\frac{p}{q}\pi) +\sin(\frac{p }{q}\pi)  \sin \frac{\pi}{2}=0 \\
 B(\varphi_0,\beta_0;\delta_0)  &=&   \sin(-\frac{p}{q}\pi) +\sin(\frac{p }{q}\pi)  \sin (\frac{\pi}{2} +2(m+1)\frac{p}{q}\pi )= 0
\end{eqnarray*}
when $(m+1) p/q$ is an integer.

 Using \ref{eqn:A} we eliminate the variable $\theta$ to rewrite \ref{eqn:B} as  
\begin{equation}
\sin \varphi   - \sin \left(  \varphi +   2 (m+1) \arcsin (\delta \sin \varphi) \right ) = 0
\label{eqn:B=0}
\end{equation}
Defining $\varphi = \pi/2 + \Delta \varphi$ and $\delta = \delta_0 + \Delta \delta$ we rewrite the above equation as
\begin{equation}
\cos (\Delta \varphi)   - \cos \left( \Delta \varphi +   2 (m+1) \arcsin ((\delta_0 + \Delta\delta) \cos \Delta \varphi) \right ) = 0
\label{eqn:BB=0}
\end{equation}

For small  $0 \sim \Delta \varphi  \gg \Delta \delta $ and keeping only lower order terms we have
$$
\arcsin \left(  (\sin  \frac{p}{q} \pi + \Delta \delta)  (\cos \Delta \varphi) \right ) \sim
 \frac{p}{q} \pi 
+ \frac{1}{\cos \frac{p}{q} \pi} \Delta \delta
- \frac{1}{2} \, \frac{\sin \frac{p}{q} \pi }{\cos \frac{p}{q} \pi} (\Delta \varphi)^2
+\ldots
$$
Using the above approximation and the fact that $2(m+1) \frac{p}{q}$ is an even integer, \ref{eqn:BB=0} can be written as
$$
\cos (\Delta \varphi)  -
 \cos \left( \Delta \varphi +  \frac{m+1}{\cos \frac{p}{q} \pi} 
\left(
 2 \Delta \delta
-  {\sin \left( \frac{p}{q} \pi \right)}  (\Delta \varphi)^2
+\ldots
\right)
 \right )  = 0
$$
Now, as  
$ \displaystyle \cos(a) - \cos(a+b) =
2 \sin\left( a+\frac{b}{2} \right) \sin\left( \frac{b}{2} \right) 
$, and keeping track of the higher order terms,  we write \ref{eqn:BB=0} as
\begin{equation} 
 \frac{m+1}{\cos \frac{p}{q} \pi}
\left( \Delta \varphi +  \frac{m+1}{ \cos \frac{p}{q} \pi}  \, \Delta \delta \right ) 
\,
\left(
2 \Delta \delta
-  {\sin \left( \frac{p}{q} {\pi}  \right) } {\Delta \varphi)^2}
\right)
+\ldots = 0
\label{eqn:cubica}
\end{equation}
This expression explicitly shows the cubic bifurcation, as we have 3 solutions if $\Delta \delta > 0$ and only one if  $\Delta \delta \le 0$.
 
So, given $\delta$ close to $\delta_0$ we have a normal point given by 
$\varphi  \sim \varphi_0 - \frac{m+1}{\cos \beta_0} (\delta - \delta_0)$ 
and, for $\delta > \delta_0$, we have two other normal points  given by
$\varphi \sim \varphi_0 \pm \frac{2}{\delta_0} (\delta - \delta_0)$. 
More precisely, for each solution $\varphi$, the normal point is given by 
$\theta = \arcsin (\delta \sin \varphi)$ and $s= \varphi-\theta$. 
Moreover, if $\delta \ne \delta_0$ these normal points are transverse.
This shows that $(\omega_0,\beta_0) = (-\frac{\pi}{2},0) \in L^0 \cap G^{-1}(L^0) \subset M_{\inn}$ is a cubic tangent normal point for $\delta=\delta_0$ that unfolds generically in this parameter as illustrated in Figure~\ref{fig:tangencias}.
\end{proof}

We will look at the set $\widetilde{S}_{\delta_0,r}$ which is the closure of the  connected component of $M_{\inn}  \backslash \mathcal{S}_{\inn}^-$ 
containing  the tangent normal point $(\omega_0,\beta_0) = (- \frac{\pi}{2}, 0)$.
We observe that the definition of  $\widetilde{S}_{\delta_0,r}$ is the same as the sets $S_{\delta,r}$ 
{introduced} in the previous section for transverse normal points in $X_{\delta}$ and therefore they share some properties. 

From  Lemma \ref{lem:BTP} above and its proof, it is clear that if $\delta$ is close to $\delta_0$ the (transverse) normal points appearing in the bifurcation process  (one for $\delta < \delta_0$ or three for $\delta > \delta_0$) are close to  $(- \frac{\pi}{2}, 0)$, the tangent point for $\delta_0$.
Since the boundaries of the connected components of  $M_{\inn}  \backslash \mathcal{S}_{\inn}^-$ vary continuously with $\delta$ and $r$,  
we can adjust the set of parameters $R_{\delta_0}$ by 
choosing $\epsilon_0$ and $r_0$ small enough such that for any 
$(\delta,r) \in R_{\delta_0}$ these normal points are in the connected component 
of $ M_{\inn}\backslash\mathcal{S}_{\inn}^-$ containing
the point
$ (- \frac{\pi}{2}, 0)$.
However, due to constructions that will  intervene later, we will eventually need to take smaller $\epsilon_0$ and $r_0$.
These components are denoted by
$\widetilde{S}_{\delta,r}$ 
and their images  
$\widetilde{U}_{\delta,r} =   G (\widetilde{S}_{\delta,r}) \subset M_{\inn}  \backslash \mathcal{S}_{\inn}^+$. 
It is worthwhile to remember that all points in the same connected component have the same returning time characterized by $m$.
 
In what follows, we describe the set  $\widetilde{S}_{\delta_0,r}$ for an arbitrary $r \le r_0$
(as usual we will drop the subscripts in sets and maps when the identification is obvious). 
We stress that $r_0$ is to be chosen small enough in order that all the arguments and the description bellow applies even for different values of $\delta$.

\begin{figure}[h]
\includegraphics[width=0.25\hsize]{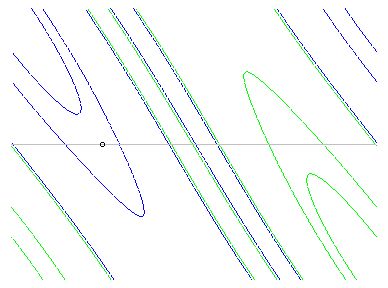}
\hfill
\includegraphics[width=0.25\hsize]{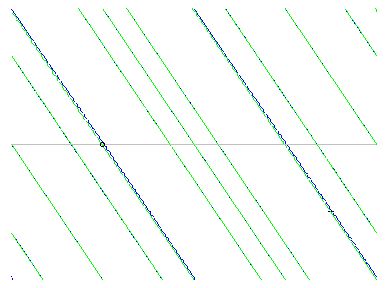}
\hfill
\includegraphics[width=0.25\hsize]{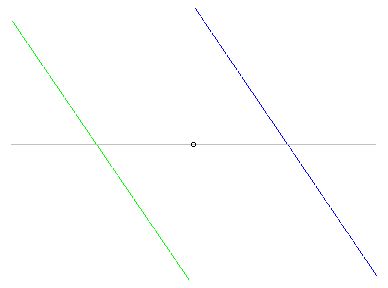}
\vskip 0.3cm
\includegraphics[width=0.25\hsize]{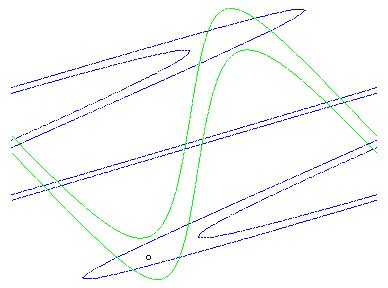}
\hfill
\includegraphics[width=0.25\hsize]{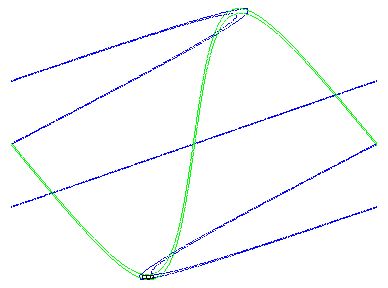}
\hfill
\includegraphics[width=0.25\hsize]{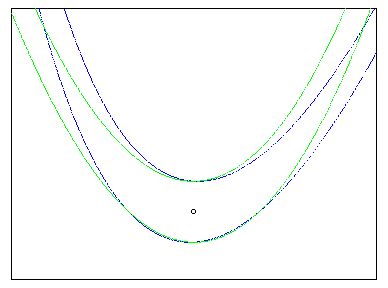}
\caption{The sets $\widetilde S$ (top)  and  $T(\widetilde S)$ (bottom) for $r<delta$ and $r \approx 0$ (also in zoom)} 
\label{fig:Stilde}
\end{figure}

Since the initial observations in the proof of Lemma \ref{lem:strips} do not rely on transversality 
of the normal point, they also apply here.
For $\delta = \delta_0$, the point  $(\omega_0,\beta_0) = (-\frac{\pi}{2},0) \in \widetilde{S}_{\delta_0,r}$ is a tangent normal point with return to the obstacle time $m+2$. 
As noticed in the proof of Lemma $\ref{lem:strips}$, 
for $r$ small enough,
the image  $T(\widetilde{S}_{\delta_0,r}) \subset M_{\out}$ 
is  the connected component of $M_{\inn}^+ \cap F^{-m}(M_{\inn}^-)$
containing the point $(\omega_0,\beta_0) = T(-\frac{\pi}{2}, 0) \in L_{\delta_0}^+ \cap F^{-m}(L_{\delta_0}^-)$.  
Moreover,  the curves in $\partial M_{\inn}^+$ are 
 ${\mathcal C}^1$ close to $L_{\delta_0}^+$ while the curves in $F^{-m}(\partial M_{\inn}^-)$ are ${\mathcal C}^1$ close to the curve $F^{-m}(L_{\delta_0}^-)$. 
As $L_{\delta_0}^+$ and  $F^{-m}(L_{\delta_0}^-)$ are in fact 
topologically transverse,
for small $r$, the boundary of $T(\widetilde{S}_{\delta_0,r})$ contains two curves
belonging to different components of $F^{-m}(\partial M_{\inn}^-)$ with endpoints in $\partial M_{\inn}^+$. 
Hence $\widetilde{S}_{\delta_0,r}$ is  a strip bounded by two curves in the singular set 
$G^{-1}_{\delta_0,r} ( \partial M_{\inn})$ connecting the two components of $\partial M_{\inn}$. 
This can be observed in Figure \ref{fig:Stilde}.

The description above implies that $T(\widetilde S_{\delta_0,r})$ converges to the point 
$T(-\frac{\pi}{2},0) \in M_{\inn}^{+}\backslash H_{\delta_0}$ as $r \to 0$.
So, for $r$ small enough the set $T(\widetilde{S}_{\delta_0,r})$
itself is contained in $M_{\inn}^{+}\backslash H_{\delta_0}$, which implies that 
$\widetilde{S}_{\delta_0,r} \subset M_{\inn}\backslash H^{-}_{\delta_0,r}$.

Even thought $\widetilde{S}_{\delta_0,r}$ is a strip, it {may not} be essentially a parallelogram as its boundaries {may not} be stable curves. 
However, the connected component of $M_{\inn} \backslash H^-_{\delta_0,r}$ containing  $\widetilde{S}_{\delta_0,r}$ is a stable strip 
when $r \sim  0$ since the two curves of $\partial H^-_{\delta_0,r}$ connecting   $\partial M_{\inn}$ are uniformly $C^1$ close to the lines $|\sin(\omega+\beta)|= \delta_0$.

Analogous properties can be derived for the set 
$\widetilde{U}_{\delta_0,r}  \ni (\frac{\pi}{2},0) = G(-\frac{\pi}{2},0)$:
For small $r$, $\widetilde{U}_{\delta_0,r}  \subset  M_{\inn} \backslash H^+_{\delta_0,r}$ is a strip bounded
by two curves in different components of $G(\partial M_{\inn})$ connecting $\partial M_{\inn}$.  Although $\widetilde{U}_{\delta_0,r}$ is not an unstable strip, it is contained in the unstable strip $M_{\inn} \backslash H^+_{\delta_0,r}$.

As we have a  continuous dependence 
of maps and sets
on $(\delta,r) \in R_{\delta_0}$ ,  the properties described above for  $\delta=\delta_0$  hold  
for all sets $\widetilde{S}_{\delta,r}$ and  $\widetilde{U}_{\delta,r}$  for any $(\delta,r) \in R_{\delta_0}$ as long as
$\epsilon_0$ and $r_0$ are properly chosen.


The geometric conditions stated bellow will provide the technical tools to prove the existence of homoclinic tangencies, 
as they ultimately will relate the behavior of segments of the stable manifold to the curve $G^{-1}(L^0)$ in the neighborhood $\widetilde S$ of the tangent normal point.

\begin{lemma} We can choose $\epsilon_0$ and $r_0$ such that for any $(\delta,r) \in R_{\delta_0}$ the local stable manifold of any point in the corresponding hyperbolic set 
$\Lambda_{\delta,r}$  has a component  connecting the two curves in $\partial \widetilde{U}_{\delta,r} \cap G_{\delta,r}( \partial M_{\inn})$. 
Moreover this component does not intersect the line $L^0$.
\label{lem:geometrico}
\end{lemma}

\begin{proof}
From  Section \ref{sec:hyperbolic}, if  $(\delta,r) \in R_{\delta_0}$,  the 
local stable manifold of any point in $\Lambda_{\delta,r}$  is a stable curve  inside $H_{\delta,r}^-$ and connecting the two components of $\partial M_{\inn}$.
Moreover, for small $r$,  the boundary $\partial H_{\delta,r}^-$ is close to the straight decreasing lines  
$ \vert \sin(\omega+\beta) \vert = \delta $ and so the local stable manifolds are inside the region 
$ \vert \sin(\omega+\beta) \vert < \delta $ 

It is clear that for $\delta$  close to  $\delta_0$ and small $r$, the boundary $\partial H_{\delta,r}^-$  belongs to a small tubular neighborhood of the lines  $\vert \sin(\omega+\beta) \vert = \delta_0 $. 
Thus we can take smaller $\epsilon_0$ and $r_0$ 
so that,  for any $(\delta,r) \in R_{\delta_0}$, the  set $\Lambda_{\delta,r}$ and its {\em local stable foliation} are contained in the interior of the two strips defined by $|\sin(\omega+ \beta)| < \delta_0$.

On the other hand, the boundary $\partial H_{\delta,r}^+$ is close to the lines $\vert \sin(\omega-\beta) \vert = \delta_0 $
and again we can adjust  $\epsilon_0$ and $r_0$ such that the set $\tilde{U}_{\delta,r} \subset M_{\inn}\backslash H^+_{\delta,r}$ is in the interior of the narrow strip defined by 
{$ \sin(\omega-\beta)  > \delta_0$  containing $(\pi/2,0)$} for any $(\delta,r) \in R_{\delta_0}$.
We note that this strip is crossed by the two {\em decreasing} strips $|\sin(\omega+\beta)| < \delta_0$ (see Figure \ref{fig:lemageometrico}).  

This  implies that  the local manifolds must cross the strip  $\sin(\omega-\beta) \geq \delta_0$ and so the strip $\tilde{U}_{\delta,r}$ which is inside it.
In particular, 
any local stable manifold of 
$\Lambda_{\delta,r}$ must have an arc connecting the two components of $\partial \tilde{U}_{\delta,r} \cap G(\partial M_{\inn})$. 

\begin{figure}[h]
\begin{center}\includegraphics[scale=0.4]{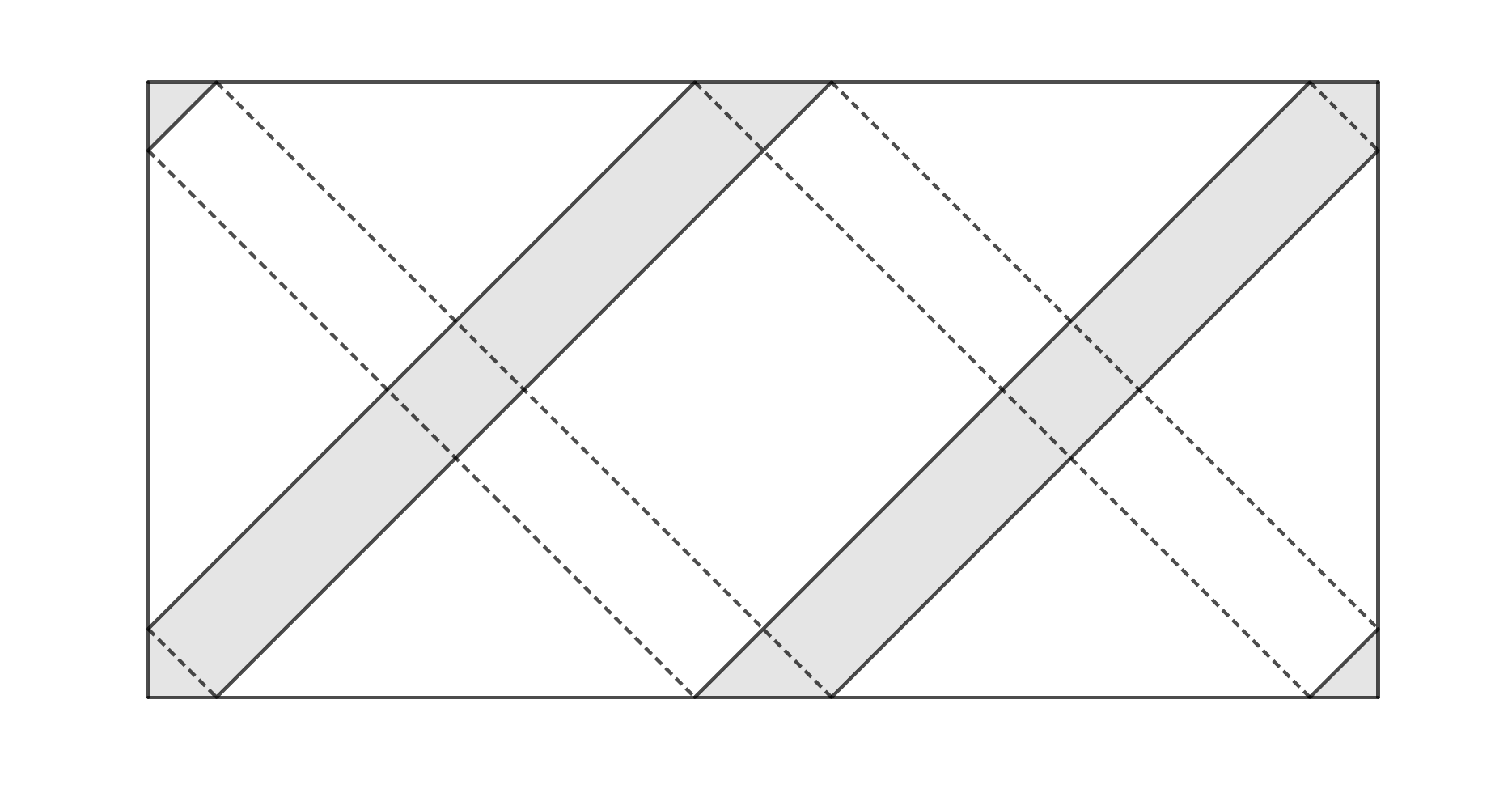} \end{center}
\caption{
Schematic representation of the geometric construction in the proof of Lemma~\ref{lem:geometrico}
($ |\sin (\omega+ \beta)|=\delta$ (dotted line),
$ |\sin (\omega- \beta)|=\delta$ (solid line),
$|\sin (\omega- \beta)| \geq \delta$ (gray region))
}
\label{fig:lemageometrico}
\end{figure}

Finally, we observe that  the set $L^0 \cap \tilde{U}_{\delta,r}$ belongs to the intersection between the strips $|\sin(\omega+\beta)| > \delta_0$ and  $|\sin(\omega-\beta)| > \delta_0$. 
As any local stable manifold 
belongs to the region $|\sin(\omega+\beta)| < \delta_0$, we have that its intersection with $\tilde{U}_{\delta,r}$ is disjoint from $ L^0$. 
\end{proof}


The strategy of the following lemma, is to obtain points of quadratic tangency between a stable manifold $W^s_{\delta,r}$ and the horizontal line $L^0$, for a set of parameters. 
The symmetry of the phase space,  implies that these points also correspond to tangencies between stable and unstable manifolds,
since the image of the stable manifold of a point by the involution is the unstable manifold  its symmetric. 
Thus  we actually have  heteroclinic tangencies. 
{Considering  the stable manifold of
symmetric periodic points produces homoclinic quadratic tangencies that unfold generically in $r$.}
By {\em unfolding a quadratic tangency}
we mean that,  if $W^s_{\delta,r^*}$ is tangent to $L^0$, for  $r<r^*$ $W^s_{\delta,r}\cap L^{0}= \emptyset$ and for $r>r^*$  $W^s_{\delta,r}\cap L^{0}$ has two distinct points (or the other way around).

\begin{lemma} \label{lem:tanquad}
Let $\delta_0 = \sin(\frac{p}{q}  \pi)$, and 
$z_{\delta,r} \in \Lambda_{\delta,r}$
be a two parameter continuous family of symmetric periodic points for $(\delta,r) \in R_{\delta_0}$.
Then, there is a curve of parameters $\Gamma \subset R_{\delta_0}$ such that if $(\delta,r) \in \Gamma$
{the stable and unstable manifolds of the point $z_{\delta,r}$ have a quadratic tangency}
which unfolds generically by fixing $\delta$ and varying $r$.
\end{lemma}

\begin{proof} 
{For  $(\delta,r) \in R_{\delta_0}$, let ${\mathcal W}_{\delta,r} = W^s_{loc}(z_{\delta,r})$ 
be the local stable manifold of the symmetric periodic point $z_{\delta,r}$.
The curves ${\mathcal W}_{\delta,r}$ connect the two components of $\partial M_{\inn}$ and converge in the ${\mathcal C}^1$ topology, as  $(\delta,r) \to (\delta_0,0)$, to the line 
$ \omega+\beta= \omega_{\delta_0}$ where $(\omega_{\delta_0},0) \in L^0$
is a transverse normal point. 
By Lemma \ref{lem:geometrico} we have that
${\mathcal W}_{\delta,r} \cap \tilde{U}_{\delta,r}$ is an arc connecting the two curves of 
 $\partial \widetilde{U}_{\delta,r} \cap G_{\delta,r}(\partial M_{\inn})$. Moreover, this arc does not intersect $L^0$.
Its inverse image $G^{-1}_{\delta,r}({\mathcal W}_{\delta,r}) $ 
is a curve  in $\widetilde{S}_{\delta,r}$ connecting $\partial M_{\inn}$
and so intersects  $L^0$. 
We will show that, for a curve $\Gamma$ of parameters in $R_{\delta_0}$,  this intersection  
is a quadratic homoclinic tangency unfolding generically in the parameter $r$ (as usual we are dropping some subscripts).
}

The idea of the proof is to repeat the construction of Lemma \ref{lem:BTP} including the effect of varying $r$ in the neighborhood of the cubic tangency. 
In order to obtain the unfolding of the tangency
we observe,  as in the proof of Lemma \ref{lem:BTP}, 
that points in $G^{-1}({\mathcal W}_{\delta,r}) \cap \widetilde{S}_{\delta,r} \cap L^0$ correspond to the intersection
\begin{equation}
T (G^{-1}({\mathcal W}_{\delta,r}) \cap L^0) = F^{-m}\circ T^{-1} ({\mathcal W}_{\delta,r}) \cap T(L^0)=  F^{-m}\circ T^{-1} ({\mathcal W}_{\delta,r}) \cap L_{\delta}^+
\label{eqn:I1}
\end{equation}
As for $(\delta,r)$ close to $(\delta_0,0)$, ${\mathcal W}_{\delta,r}$ is close to the line
 $\omega+\beta= \omega_{\delta_0}$, there is a smooth function $\epsilon(\beta;\delta,r)$ 
such that  ${\mathcal W}_{\delta,r}$ can be written as
$$
\omega+\beta= \omega_{\delta_0} + \epsilon(\beta;\delta, r)
\hbox{ \ \ where   \ \ } \epsilon(\beta;\delta,r) \rightarrow 0  \hbox{ as  } (\delta,r) \rightarrow (\delta_0,0)
$$
The preimage $T^{-1}({\mathcal W}_{\delta,r}) \subset M_{\out}$ is a curve connecting $\partial M_{\inn}^-$ defined by
\begin{eqnarray*} 
&& \sin\theta +\delta \sin(\theta - s)= - r\sin\beta 
\\
&& 2\beta = \theta - s + \omega_{\delta_0}+ \epsilon(\beta;\delta,r)
\end{eqnarray*} 
and 
 $F^{-m}\circ T^{-1}({\mathcal W}_{\delta,r}) \subset M_{\out}$ is written as 
$$ \sin\theta +\delta \sin(\varphi-2(m+1)\theta) = - r \,D(\varphi,\theta;\delta,r)
$$
where $\varphi = s + \theta$ and 
$$ D(\varphi,\theta;\delta,r)
=\sin \left(\dfrac{\omega_{\delta_0}+ \epsilon(\beta;\delta,r) + \varphi -2(m +1)\theta + m\pi }{2}  \right ) 
$$
So the intersection $L_{\delta}^+ \cap F^{-m}\circ T^{-1}({\mathcal W}_{\delta,r})$ is {a  solution} of the following system,  whose left hand side is the same of  the one considered in Lemma $\ref{lem:BTP}$ (equations \ref{eqn:A} and \ref{eqn:B}) .
\begin{eqnarray}
A(\varphi,\theta;\delta) = \sin\theta +\delta \sin\varphi &=&0  \nonumber\\
B(\varphi,\theta;\delta) = \sin\theta +\delta \sin(\varphi-2(m+1)\theta) &=& - r \, D(\varphi,\theta;\delta,r) 
\label{eqn:S2}
\end{eqnarray}
At lower order, \ref{eqn:S2} is equivalent to the following cubic equation, which should be compared to \ref{eqn:cubica}.
\begin{equation} 
 \frac{m+1}{\cos \frac{p}{q} \pi}
\left( \Delta \varphi +  \frac{m+1}{ \cos \frac{p}{q} \pi}  \, \Delta \delta \right ) 
\,
\left(
2 \Delta \delta
-  {\sin \frac{p}{q} \pi } {\Delta \varphi^2}
\right)
+\ldots = - r D_0
\label{eqn:cubica2}
\end{equation}
where, using that $\epsilon(0;\delta_0,0)=0$, 
$$D_0=D\left(\frac{\pi}{2}, -\frac{p}{q}\pi, \delta_0,0 \right)
=\sin \left( \frac{\omega_{\delta_0}+  \frac{\pi}{2}}{2}- (m+1)\frac{p}{q}\pi + m\frac{\pi}{2}\right) 
$$

{
It is important to notice that, according to Lemma \ref{lem:geometrico}, the function $D$ cannot be 0 in the neighborhood considered here,  since $D=0$ would imply the existence of a point  $L^0 \cap {\mathcal W_{\delta,r}} $ in $\widetilde U_{\delta,r}$.  In particular, $D_0 \ne 0$.}

We can conclude that for $\delta$ close to $\delta_0$ and $r$ small, the curves {
$G^{-1}({\mathcal W}_{\delta,r})$, in the neighborhood of $(-\frac{\pi}{2},0)$, are essentially  translations of the curve $G^{-1}(L^0)$}

 For each $\delta $ we can adjust this translation to produce the unfolding of a quadratic tangency between the stable manifold and the horizontal symmetry line, 
as in Figure~\ref{fig:quadratic}.

\begin{figure}[h]
\includegraphics[width=0.24\hsize]{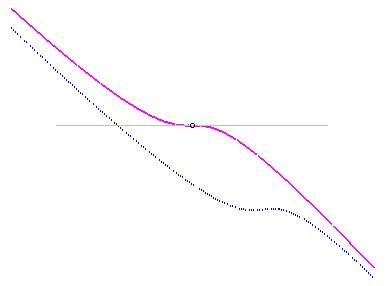}\hfill
\includegraphics[width=0.24\hsize]{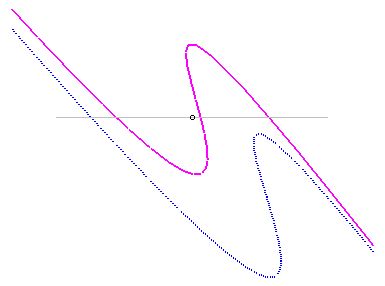}\hfill
\includegraphics[width=0.24\hsize]{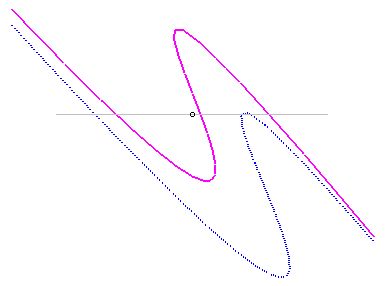}\hfill
\includegraphics[width=0.24\hsize]{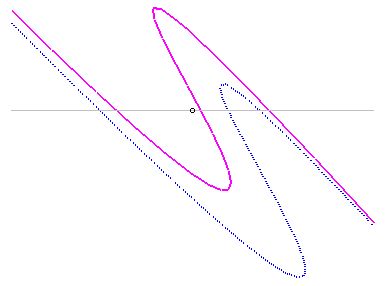}
\caption{Bifurcation of Homoclinic Tangency. ($L^0$, $G^{-1}(L^0) $ and $W^s$ for $\delta=\delta_0$ and
$\delta > \delta_0$ with three different values of $r$).}
\label{fig:quadratic}
\end{figure}

In fact, a quadratic tangency occurs if a solution of \ref{eqn:cubica2} also satisfies
$$
2 \Delta \delta
-  \delta_0 {\Delta \varphi^2} - 2 \delta_0\,  \varphi \, (\Delta \varphi +  \frac{m+1}{\cos \frac{p}{q} \pi}  \Delta \delta)  = 0
$$
Solving this equation and  \ref{eqn:cubica2} for $ \Delta \delta$ and $r$, we obtain,  at lower order  
\begin{eqnarray}
 \delta &\approx& \delta_0 + \frac{3}{2} \, {\delta_0 \, (\varphi - \frac{\pi}{2})^2}  \label{eqn:gamma} \\ 
r &\approx& 2\, \frac{(m+1) \tan \frac{p}{q}\pi}{- D_0} \, (\varphi-\frac{\pi}{2})^3 \nonumber
\end{eqnarray}
The two equations above define a curve $\Gamma$ in the parameter set,  
approaching the point $(\delta_{0},0)$ from $\delta > \delta_0$, along which, the curve 
{$G^{-1}({\mathcal W}_{\delta,r})$} and the line $L^0$ have a quadratic tangency unfolding generically with $r$.
As that ${\mathcal W}_{\delta,r}$ is constructed from the stable manifold of symmetric periodic points, these tangencies are in fact homoclinic. 
\end{proof}

At this point it is worthwhile to notice that families constructed from different symmetric periodic points,
will give rise to different homoclinic tangency curves
also abutting $(\delta_0,0)$.

{
Before proceeding to the proof of Theorems \ref{thmB} and \ref{thmC} we will show that symmetric periodic points indeed exist in the hyperbolic set. 
To this end, we refer to the construction of the set $\Lambda_{\delta,r}$ with fixed parameters $(\delta,r) \in \Omega_0$.
This construction is based on the strips  $S_i$ and $U_i$ for $i=1,\ldots,n_\delta$ associated to the transverse normal points in $X_{\delta}$ (Lemma \ref{lem:extendedhorseshoe}).
where a point in $\Lambda_{\delta,r}$ is specified by its symbolic representation, which is a sequence 
$a \in \tilde \Sigma \subset \Sigma = \{ 1, \ldots, n_{\delta} \}^{\mathbf Z}$. 
Let us consider a set of strips  $S_{a_0},S_{a_2}, \ldots S_{a_k}$ with $a_i \in \{ 1, \ldots, n_{\delta} \}$ and such that  $U_{a_i} = G( S_{a_i}) $ crosses $S_{a_{i+1}}$.
This equivalent to the say that the word $[a_0....a_{k}]$
appears in some  sequence
$a \in \tilde \Sigma$.


The horizontal curve $ U_{a_{k}}\cap L^0$  connects the two components of $G(\partial M_{\inn}) \cap U_{a_{k}}$. 
It follows  from the properties of stable and unstable strips that the preimage $\lambda_{k}= G^{-1}(U_{a_{k}}\cap L^0) \subset S_{a_{k}}$ is a stable curve connecting the two components of $\partial M_{\inn} \cap S_{a_{k}}$. Thus $\lambda_{k} \cap U_{a_{k-1}}$ is a stable curve connecting the two components of $G(\partial M_{\inn}) \cap U_{a_{k-1}}$ and hence $\lambda_{k-1}= G^{-1}(\lambda_{k} \cap U_{a_{k-1}})$ is a stable curve in $S_{a_{k-1}}$ connecting $\partial M_{\inn}$. 
Iterating this construction we define, for $j=0,...,k$, the stable curves $\lambda_{k-j} \subset S_{a_{k-j}} \cap G^{-j}(L^0)$  each of which connects  $\partial M_{\inn}$.
Thus the the unstable curve $\lambda_0 \subset \bigcap_{j=0}^k G^{-j}(S_{a_{j}})$ intersects transversally the horizontal curve $L^0 \cap S_{a_{0}}$.  Since $\lambda_0 \in G^{-k-1}(L^0)$, the intersection $z=\lambda_{0} \cap L^0 \subset S_{a_0}$ is a symmetric periodic point 
having $[a_0...a_{k}]$ in its symbolic representation. 

It is clear that this construction produces a  continuous two parameter family  of symmetric periodic points 
$z_{\delta,r} \in L^0 \cap G_{\delta,r}^{-k-1}(L_0)$ for $(\delta,r) \in R_{\delta^*}$ for some $\delta^*$.



{Now, Lemma \ref{lem:tanquad} immediately gives Theorem \ref{thmB}}.

{\bf Theorem \ref{thmB}:}
{\em There is a set $\Omega'_0 \subset \Omega_0$ 
accumulating the point $(1,0)$ such that if $(\delta, r) \in \Omega_0'$, then  $G_{\delta,r}$ presents  quadratic homoclinic tangencies unfolding generically with the parameter $r$}


\begin{proof}[Proof of Theorem $\ref{thmB}$]
To each $ \delta_{\kappa}=\sin \kappa\pi $,  with  $\kappa$ a rational number in $(\frac{1}{4},\frac{1}{2})$, 
there is a set $R_{\kappa} = (\delta_{\kappa}-\epsilon_{\kappa}, \delta_{\kappa} + \epsilon_{\kappa} ) \times  (0, r_{\kappa}] \in  \Omega_0$ 
where lemmas \ref{lem:extendedhorseshoe}, \ref{lem:geometrico} and \ref{lem:tanquad} 
hold.  

For a fixed arbitrary $\delta_{\kappa}$ 
as above, we can pick a {family of symmetric periodic points}
 in the hyperbolic set  $\Lambda_{\delta_{\kappa},r}$ to obtain, by Lemma \ref{lem:tanquad}, a
curve $\Gamma \subset R_{\kappa}$ of homoclinic tangencies. 

We define $R_{\kappa}' \subset R_{\kappa}$ as the the union of the curves  such that $G_{\delta,r}$ with $(\delta,r) \in R_{\kappa}'$ unfolds generically a quadratic homoclinic tangency.
{It is worth mentioning that, from  \ref{eqn:gamma}, $R_{\kappa}'$ is contained in $\delta > \delta_\kappa$ and abutts $(\delta_{\kappa},0)$.}

The set 
 $ \left\{
 \delta_{\kappa}=\sin(\kappa\pi) \hbox{  , with } \kappa \in (\frac{1}{4},\frac{1}{2}) \cap \mathbb{Q}
\right \}$ is dense in $(\frac{1}{\sqrt{2}} ,1)$ and
the set  $\Omega_0'= \bigcup_{\kappa} R_{\kappa}'$ accumulates the point $(1,0)$. This concludes the proof of the theorem.
\end{proof}

\begin{remark}
{To each continuous two parameter family of symmetric periodic points in $ \Lambda_{\delta,r}$ it corresponds a tangency curve $\Gamma$. 
The union of these curves in each $R_k$ is a set of tangency bifurcations with an {intricate geometric structure} that we  do not intend to describe here.
}
\end{remark}


We close this section with the proof of our third theorem.  

\begin{proof}[Proof of Theorem $\ref{thmC}$]  
{Consider a set of parameters $R_{\kappa}$, as in the  in the proof of Theorem \ref{thmB}.
Let $E_{\kappa}$ be the subset of all pairs $(i,j) \in \{1,...,n_{\delta_{\kappa}}\}^2$ such that $U_{i}=G(S_{i})$ crosses $S_{j}$.
For $(\delta,r) \in  R_{\kappa}$, as in Section \ref{sec:hyperbolic},  the strips  
$S_1,...S_{ n_{\delta_{\kappa}}}$ are associated to the hyperbolic set. 
In particular, $\Lambda_{\delta,r} \subset \bigcup_{(i,j) \in E_{\kappa}} U_{i} \cap S_{j} $.
}

The existence of elliptic periodic points follows from a homoclinic bifurcation
associated to a continuous family of 
{specifically chosen} symmetric periodic points $y_{\delta,r} \in \Lambda_{\delta,r}$.
The points $y_{\delta,r}$ are constructed from a given admissible word $[a_0....a_{m}] \in \{1,...,n_{\delta_{\kappa}}\}^k$,  as explained {earlier}. 

We can choose a word $[a_0...a_{m}]$ containing {every}
admissible sequence of two symbols of the form $[a_{i}a_{j}]$.
The orbit  of the resulting point $y_{\delta,r}$ visits {all the components}  $U_{i}\cap S_{j}$ with
 $(i,j) \in E_{\kappa}$, spreading over the hyperbolic set.
This also implies that, as  $(\delta,r) \rightarrow (1,0)$, the orbit of $y_{\delta,r}$ tends to fill the entire phase space. 
More precisely,
the maximum distance of points of phase space to the union point of the orbit  of $y_{\delta,r}$ goes to $0$ as  $(\delta,r) \rightarrow (1,0)$

Fixing $\delta^* \in (\delta_\kappa,\delta_\kappa + \epsilon_\kappa)$ we consider the one parameter family 
of maps $G_{\delta^*,r}$ and the related family of symmetric periodic points $y_{\delta^*,r}$  with $r \in (0, r_{\kappa}]$. 
From Lemma \ref{lem:tanquad} there is  $r^*$ such that the invariant manifolds of $y_{\delta^*,r^*}$ have a quadratic homoclinic tangency  unfolding generically  in the parameter $r$. 
From Duarte's result \cite{duarte} there is a subset $ I \subset (0,r_{\kappa}]$ accumulating $r^*$, such that for every $r \in I$, 
the closure of the  generic elliptic periodic points of $G_{\delta^*,r}$ contains the orbit of  $y_{\delta^*,r}$. 

Thus, in each $R_{\kappa}$ there is {a subset  of parameters} $R''_{\kappa}$ for which the map $G_{\delta,r}$ has a set $\mathcal{E}_{\delta,r}$ of 
{generic elliptic periodic points}.
Clearly $R''_{\kappa}$ accumulates the set $R_{\kappa}'$ {where homoclinic tangencies do exist}.  

{Finally}, the set of parameters $\Omega_0''= \bigcup_{\kappa} R_{\kappa}'' \subset \Omega_0$  accumulates $(1,0)$.
For each $(\delta,r) \in \Omega_0''$, the set of generic elliptic points $\mathcal{E}_{\delta,r}$  accumulates  the orbit of symmetric periodic points $y_{\delta,r}$. 
{This fact, together with the properties of the orbit of $y_{\delta,r}$, implies} that the maximum distance of points in the phase space to the set 
$\mathcal{E}_{\delta,r}$ also goes to zero as $(\delta,r) \rightarrow (1,0)$.
\end{proof}